\documentclass[a4paper]{article}

\usepackage{graphicx}
\usepackage{amscd}
\usepackage{amsmath}
\usepackage{amsfonts}
\usepackage{amssymb}
\usepackage{subcaption}
\usepackage[active]{srcltx}
\usepackage{theorem}
\usepackage{enumerate}
\usepackage{a4wide}
\theorembodyfont{\rmfamily}

\newcommand{\R}{I\!\!R}
\newcommand{\N}{\mathbb N}
\newcommand{\C}{\mathbb C}

\newtheorem{lem}{Lemma}[section]

\newtheorem{pro}[lem]{Proposition}

\def\proof{\par{\bf Proof:} \ignorespaces}
\def\endproof{{\ \vbox{\hrule\hbox{%
   \vrule height1.3ex\hskip0.8ex\vrule}\hrule

    }}\par}

\newcommand{\yd}{y^{\delta}}

\newcommand{\tA}{{\tilde A}}

\begin{document}

\title{An entropic Landweber method for linear ill-posed problems}

\author{Martin Burger, Elena Resmerita and Martin Benning}

\maketitle

\begin{abstract}
The aim of this paper is to investigate the use of a Landweber-type method involving the Shannon entropy for the regularization
of linear ill-posed problems. We derive a closed form solution for the iterates and analyze their convergence behaviour both in a case of reconstructing general nonnegative unknowns as well as for the sake of recovering probability distributions.  Moreover, we discuss several variants of the algorithm and relations to other methods in the literature. The effectiveness of the approach is studied numerically in several examples.
\end{abstract}

\section{Introduction}

This work deals with   linear ill-posed equations $Au=y$ with $A:X\rightarrow Y$ acting between a Banach space $X$ and a Hilbert space $Y$, for which solutions with specific properties (such as positivity) are sought.  In this respect, we consider   iterative regularization methods of the following type
 \begin{equation}\label{method0}
u_{k+1}\in \mbox{arg}\min_u\left\{\frac{1}{2}\|Au-y\|^2+cd(u,u_{k})-\frac{1}{2}\|Au-Au_{k}\|^2\right\},\,\,\, k\in\N,
\end{equation}
where $d=D_f$ denotes the Bregman distance \cite{Bregman} associated with a convex functional $f:X\rightarrow\mathbb{R}\cup\{+\infty\}$ which is nonnegative and $c$ is some positive number. The term $d(u,u_{k})$ acts as a penalty enforcing the desired features for the solutions. 

\noindent Note that we can rewrite \eqref{method0} as
 \begin{equation}\label{method0a}
u_{k+1}\in \mbox{arg}\min_u\left\{\langle Au-y,Au_{k}-y\rangle +cd(u,u_{k}) \right\},
\end{equation}
which shows that the scheme can be obtained also by linearizing the quadratic data-fitting term $\displaystyle{\frac{1}{2}\|Au-y\|^2}$ at the current iterate $u_k$. 
This class of methods  incorporates several procedures that have been proposed so far in the literature. For instance, the classical  case when $f$ is  quadratic in Hilbert spaces reduces to the Landweber method, as emphasized by \cite{DDD} and as investigated for nonlinear operator equations by means of surrogate functionals in \cite{RamTes},  see also the discussion in \cite{Scherzer}. The case of quadratic $f$ in reflexive Banach spaces has been studied by  \cite{SchSchLou}.  The setting when $f$ is the total variation functional smoothed by a quadratic has been analyzed by \cite{BacBur}, requiring fine analysis tools due to the bounded variation function space context. The case of $\ell^1$-penalties has been treated in \cite{Yin_etal, Cai_etal},  resulting in the so-called  linearized Bregman algorithm. In all those cases however, some quadratic term had to be part of $f$ to guarantee even well-definedness of the iterates and subsequently convergence. 

\noindent We are interested here in the Shannon entropy setting without any quadratic term, i.e.
$$ f(u) = \int_\Omega u(t) \ln u(t) ~dt. $$
We mention that one can alternatively consider a linear shift to 
$$ \tilde f(u) = \int_\Omega u(t) \ln u(t) - u(t) +1 ~dt, $$
which is a nonnegative functional inducing the same Bregman distance. 

\noindent This raises challenges to analyze the problem in the $L^1$ setting without quadratic terms in the functional, but provides a simple closed iterative method with preserving the sign of the starting point (function) along the iterations. The latter formulation involves entropic projections, as one can see in the following section.

Moreover, we shall also be interested in the solution of inverse problems with unknowns being probability densities, i.e. we minimize on the domain of $f$ subject to the constraint
\begin{equation}
	\int_\Omega u(t)~dt = 1, 
\end{equation}
which again results in a simple closed iterative form. 

The advantages of using  Bregman projections for solving variational problems with unknown probability densities have been exploited by several authors before, e.g. in optimal transport (cf. \cite{benamou,peyre}). 

\noindent In order to write both problems in a closed form, we will use the equivalent formulation
 \begin{equation}\label{method0b}
u_{k+1}\in \mbox{arg}\min_{u}\left\{\langle Au-y,Au_{k}-y\rangle +cd(u,u_{k}) +\chi_m(u)\right\},\,\,\,m\in\{0,1\},
\end{equation}
where $\chi_0 \equiv 0$ denotes the original problem without integral constraint, and 
$$ \chi_1(u) = \left\{  \begin{array}{ll} 0 & \mbox{if~} \int_\Omega u(t)~dt  =1,\\ + \infty &\mbox{else},\end{array}\right. $$
is employed for enforcing probability densities. The minimization is taken here over the domain of the entropy functional.

One finds the above entropy based algorithm in the finite dimensional optimization literature, as well as in the machine learning one, under quite different names. One could mention the mirror descent type algorithms for function minimization introduced in \cite{NemYud} and   the Bregman-distance version with emphasis on entropy  in \cite{BecTeb}, and the exponentiated gradient descent method for linear predictions - see \cite{KivWar}. The work \cite{Ius} investigated three versions of the so-called approximate (linearized) proximal point methods for optimization in combination with line search strategies. The reader is referred to Sections 6.6 - 6.9 in \cite{cen_zen} for other iterative optimization methods employing the Shannon entropy.

The main contribution  of our work is  the convergence of the iterates \eqref{method0b} to a solution of the equation $Au=y$ even in an infinite dimensional setting  of such a nonquadratic penalty version, by stating also error estimates in the sense of a distance between the solution and the iterates, as opposed to the classical situation encountered in optimization, where the error for the objective function values is highlighted. 

This manuscript is organized as follows. Section 2 provides the necessary background on entropy functionals, as well as on well-definedness of the proposed iterative procedure. Section 3  analyzes (weak) convergence of the method when both a priori and a posteriori stopping rules are considered, while Section 4 deals with error estimates only for the former rule.  Section 5 explores a version of the entropic Landweber method for nonquadratic data fidelity terms. The theoretical results are  tested in Section 6 on several integral equation examples, in comparison with the Expectation-Maximization algorithm  and the projected Landweber method  - see \cite{cla_kal_res} for an overview on regularization methods for nonnegative solutions of ill-posed equations.

\section{Preliminaries}

In the following we collect some basic results and assumptions needed for the analysis below. We start with properties of the entropy and then proceed to the operator $A$.

\subsection{Entropy and Entropic Projection}

Let $\Omega $ be an open and bounded subset of ${\mathbb{R}}^d$.
\textit{The   negative of the Boltzmann-Shannon entropy}  is  the function
$f:L^1(\Omega)\rightarrow (-\infty,+\infty]$,
 given by\footnote{We use the convention $0\ln 0=0$.}
\begin{equation}\label{entropy}
f(u)=\left\{
\begin{array}{ll}
 {\int_{\Omega} u(t)\ln u(t)\,dt} & \mbox{if $u\geq 0$ a.e. and $u\ln u\in {{L}}^1(\Omega)$}, \\
&  \\
 {+\infty} & \mbox{otherwise.}
\end{array}
\right.
\end{equation}
  Here and in what follows ${{L}}_+^p(\Omega)$, $p\in [1,\infty]$, stands for the set  $\{u\in {{L}}^p(\Omega): u(t)\geq 0\,\,\,\mbox {a.e.}\}$, while $\|\cdot\|_p$ denotes, as usual, the norm of the space $L^p(\Omega)$.

The Kullback-Leibler functional or the Bregman distance  with respect to the
Bolzmann-Shannon entropy can be defined as $d:{\mbox{dom}\, f}\times {\mbox{dom}\, f}\rightarrow [0,+\infty]$ by
\begin{equation}\label{bd}
d(v,u)=f(v)-f(u)-f^{\prime}(u,v-u)
\end{equation}
where $f^{\prime}(u,\cdot)$ is the directional derivative at $u$. Here $ \mbox{dom}\,f=\{u\in {{L}}^1(\Omega): f(u)<\infty\}$ denotes the domain of $f$. One can also write
\begin{equation}\label{kl}
d(v,u)=\int_\Omega \left [v(t)\ln\frac{v(t)}{u(t)}-v(t)+u(t) \right]dt
\end{equation}
if $d(v,u)$ is finite, as one can see below.

Some properties of the entropy functionals $f$ and $d$ are recalled below (see, e.g., \cite{Res05, ResAnd07}).\smallskip

\begin{lem}\label{entropy_propr} The function defined by (\ref{entropy}) has the following properties:\smallskip
\begin{enumerate}[(i)]
\item  The domain of the function $f$ is strictly included in ${{L}}_+^1(\Omega)$.

\item  The interior of the domain of the function $f$ is empty.

\item The set $\partial f(u)$ is nonempty if and only if $u$ belongs to $L_+^{\infty}(\Omega)$ and  is  bounded away from zero. Moreover, $\partial f(u)=\{1+\ln u\}$.

\item The directional derivative of the function $f$ is given by
\[
f^{\prime}(u,v)=\int_{\Omega} v(t)[1+\ln u(t)]\,dt,
\]
 whenever it is finite.

\item For any $u,v\in {\mbox{dom}\, f}$, one has
\begin{equation}\label{inequality}
    \|u-v\|_1^2\leq
    \left(\frac{2}{3}\|v\|_1+\frac{4}{3}\|u\|_1\right)d(v,u).
\end{equation}
\end{enumerate}
\end{lem}
 
Based on Lemma \ref{entropy_propr} (iii) we define in the following
\begin{equation}\label{partialdom}
{\mbox{dom}\,\partial f}=\{u\in L^1(\Omega):  u\,\mbox{bounded and bounded away from zero} \,\,a.e.\}.
\end{equation}

\begin{lem}\label{kl_propr} The statements below hold true:
\begin{enumerate}[(i)]
\item The function $(v,u)\mapsto d(v,u)$ is convex;

\item The function $d(\cdot,u^*)$  is  lower semicontinuous  with respect to the weak topology of
$L^1(\Omega)$, whenever $u^*\in {\mbox{dom}\, f}$;

\item For any $C>0$ and any nonnegative $u\in L^1(\Omega)$, the following sets are weakly compact in $L^1(\Omega)$:
\[
\{x\in L^1(\Omega):d(x,u)\leq C\}.
\]
\item The set $\partial d(\cdot,u^*)(u)$ is nonempty for $u^*\in {\mbox{dom}\, f}$  if and only if $u$ belongs to $L_+^{\infty}(\Omega)$ and  is  bounded away from zero. Moreover, $\partial  d(\cdot,u^*)(u)=\{\ln u-\ln u^*\}$.
\end{enumerate}
\end{lem}

\noindent Denote 
\[
\langle u,v  \rangle=\int_\Omega u(t)v(t)\,dt, 
\]
for $u,v\in L^1(\Omega)$, when the integral exists.

A key observation for obtaining well-definedness of the iterative scheme as well as an explicit form for the iterates is the following result on the entropic projection.
\begin{pro} \label{explicitprop}
Let $\ell \in L^\infty(\Omega)$ and $v \in$ dom $\partial f$. Then the problem
\begin{equation}
\langle \ell, u \rangle + d(u,v) + \chi_m(u) \rightarrow \min_{u \in {\mbox{dom}\, f} }
\end{equation}
has a unique solution in the cases $m=0$ and $m=1$, respectively, given by
\begin{equation}
	u_m = c_m v e^{-\ell}, \qquad c_j = \left\{ \begin{array}{ll} 1 & \mbox{if~} m=0, \\ 
	\frac{1}{\int_\Omega v e^{-\ell}~dt }& \mbox{if~} m=1, \end{array}\right. 
\end{equation} 
which satisfies $u_m \in$ dom $\partial f$.
\end{pro}
\proof
We simply rewrite the functional as 
\begin{align*} 
\langle \ell, u \rangle + d(u,v) + \chi_m(u) &= \int_\Omega \left[ u(t) \ln \frac{u(t)}{v(t)} -u(t) + v(t) + u(t) \ell(t)\right]~dt + \chi_m(u) \\ &=
\int_\Omega \left[ u(t) \ln \frac{u(t)}{u_m(t)} -u(t) +v(t) + u(t) \ln c_m \right]~dt + \chi_m(u)\\
&=d(u,u_m) + \ln c_m \left(  \int_\Omega u(t) ~dt - 1\right) + \chi_m(u) + C_m,
\end{align*}
where $C_m=\ln c_m-\int_\Omega u_m(t)~dt$ is a constant independent of $u$. It is straightforward to notice that 
$$ \ln c_m \left(  \int_\Omega u(t) ~dt - 1\right) + \chi_m(u) = \chi_m(u). $$
Hence, the problem is equivalent to minimizing $d(u,u_m) + \chi_m(u)$. Since both terms are nonnegative and vanish for $u=u_m$, we see that $u_m$ is indeed a minimizer in ${\mbox{dom}\,f}$. Strict convexity of $d$ implies the uniqueness and since $u_m$ is the product of $v$ with a function strictly bounded away from zero it also satisfies $u_m\in$ dom $\partial f$.
\endproof

\subsection{Forward operators and entropy}

In this paper we always assume that $A:L^1(\Omega)\rightarrow Y$ is a linear and bounded operator with $Y$ being a Hilbert space. In addition to the norm boundedness of $A$, we assume a continuity property in terms of the Bregman distance. More precisely we assume that
\begin{equation} \label{condition}
	\|Au-Av\| \leq \gamma \sqrt{d(u,v)}
\end{equation}
holds on dom$(f+\chi_m)$ in the respective cases $m=0$ or $m=1$ for some positive number $\gamma$. It is easy to see that the latter is already implied by the boundedness of $A$ in case $m=1$:

\begin{lem}\label{continuity}
Let $A$ be as above with $\Vert A \Vert$ denoting its operator norm, let $u,v \in $ dom$(f+\chi_1)$, and $v \in $dom $\partial f$. Then \eqref{condition} is satisfied with $\gamma = \sqrt{2} \Vert A \Vert$.
\end{lem}
\proof
By the boundedness of $A$ we have
$$ \Vert A u - Av \Vert^2 \leq \Vert A \Vert^2 ~\Vert u - v \Vert_{L^1(\Omega)}^2. $$
Lemma \ref{entropy_propr} (v) further implies
$$ \Vert A u - Av \Vert^2 \leq \Vert A \Vert^2~ 2 d(u,v), $$
which yields the assertion.
\endproof

\medskip

\noindent We define the nonlinear functional
\begin{equation}\label{distance}
D(u,v)=cd(u,v)-\frac{1}{2}\|Au-Av\|^2,
\end{equation} 
which will be useful for the further analysis. Note that $D(u,v)\geq 0$  for $u,v \in $ dom$(f+\chi_1)$, and $v \in $dom $\partial f$, whenever $c\geq \frac{\gamma^2}{2}$ (cf. Lemma \ref{continuity}).
In case $m=0,$ we restrict the analysis to the class of operators $A$ for which $D(u,v)\geq 0$  for $u \in $ dom$f$ and $v \in $dom $\partial f$.


\section{Convergence of the Entropic Landweber Method}\label{sec:convergence}

In the following we consider the iterative method defined by \eqref{method0b}, 
where $d$ is the Kullback-Leibler divergence given by \eqref{kl}.

Noticing that $A^*$ maps to $L^\infty(\Omega)$, we can equivalently rewrite the minimization in \eqref{method0b} in the form of Proposition \ref{explicitprop}, which implies the following result.

\begin{pro} \label{definedness}
Let $u_0 \in$ dom $\partial f$. Then there exists a unique  minimizer in \eqref{method0b}   for any $k \geq 0$, given by ($\lambda = \frac{1}c$)
\begin{equation}\label{entropic}
u_{k+1}= u_{k} c_{k}^m e^{\lambda A^*(y-Au_{k})}, \qquad c_{k}^m = \left\{ \begin{array}{ll} 1 & \mbox{if~} m=0, \\ 
	\frac{1}{\int_\Omega u_{k}  e^{\lambda A^*(y-Au_{k})}~dt }& \mbox{if~} m=1, \end{array}\right. 
\end{equation}
which further satisfies $u_{k+1} \in$ dom $\partial f$.
\end{pro}

%
%

Note that, from pointwise manipulation of \eqref{entropic} we rigorously obtain the first-order optimality condition for the variational problem in each step, i.e.,
\begin{equation} \label{optimalitycondition}
	  \ln  u_{k+1}= \ln u_{k} + \ln c_k^m+ \lambda A^*(y-Au_{k}) , 
\end{equation}
where $\ln c_k^0=0$ and $\ln c_k^1$ is to be interpreted as a Lagrange multiplier for the integral constraint. In the latter case, this constant term is orthogonal to all functions of the form $w-v$, where $\chi_1(w)=\chi_1(v) = 0$. Since most estimates below for iterates will be based on taking duality products of \eqref{optimalitycondition} with such functions, they can be carried out in the same way for $m=0$ and $m=1$.

The analysis of the above method ressembles  the one for proximal point methods, which is apparent from rewriting \eqref{method0b} as
\begin{equation}
	u_{k+1}\in \mbox{arg}\min_u\left\{\frac{1}{2}\|Au-y\|^2 + \chi_m(u) + D(u,u_{k}) \right\}.
\end{equation}
 However, the quantity $D$ is neither a metric distance nor necessarily a Bregman distance of a convex function, rather a weighted difference of Bregman distances. This and the involved Kullback-Leibler divergence in an infinite dimensional setting require thus a careful  investigation.

We choose $u_0\in\mbox{dom}\,\partial f$ such that $\xi_0:=1+\ln u_0\in \mathcal{R}(A^*)$, that is $\xi_0=\lambda A^*w_0$ for some $w_0\in Y$, and denote
\begin{equation}\label{v}
v_k=w_0+\sum_{j=0}^{k-1}(y-Au_j)
\end{equation} 
for $k\geq 1$. Then \eqref{optimalitycondition} can be expressed as $\displaystyle{\xi_k=\xi_{k-1}+\ln c_{k-1}^m+\lambda A^*(y-Au_{k-1})}$ which implies
\begin{equation}\label{ve}
\xi_k=\xi_0+\sum_{j=0}^{k-1}\ln c_{j}^m+\lambda A^*\left(\sum_{j=0}^{k-1}(y-Au_j)\right)= \sum_{j=0}^{k-1}\ln c_{j}^m+\lambda A^*v_k.
\end{equation}
We show next that the entropic Landweber method converges in the exact data case.

\begin{pro} \label{conv_exact} Let $A:L^1(\Omega)\rightarrow Y$ be a bounded linear operator which satisfies \eqref{condition}  and such that the operator equation $Au=y$ has a positive solution $z$ verifying $\chi_m(z) = 0$ if $m=1$.  Let $u_0\in \mbox{dom}\,\partial f$ be an arbitrary starting element such that $1+\ln u_0\in \mathcal{R}(A^*)$. Moreover, let $\chi_m(u_0) = 0$ if $m=1$. Then the following statements are true:
\begin{enumerate}[(i)]
\item The residual $\|Au_k-y\|$ decreases monotonically.

\item The term $D(z,u_k)$ decreases monotonically.

\item  The sequences $\{u_k\}_{k\in\mathbb{N}}$ generated by the iterative method \eqref{entropic}  converge weakly on subsequences in $L^1(\Omega)$ to solutions of the equation $Au=y$,
with $\chi_m(u) = 0$ if $m=1$.
\end{enumerate}

\proof 
We will use the proximal point method techniques in order to prove the statements, by taking care of the fact that $D$ is  a nonnegative functional satisfying $D(u,u)=0$ for any $u$ in this function's domain. 

(i) We have for all $k\in \mathbb{N}$,
\[
\frac{1}{2}\|Au_{k+1}-y\|^2+D(u_{k+1},u_k)\leq \frac{1}{2}\|Au_{k}-y\|^2,
\]
which implies that the sequence $\{\|Au_k-y\|\}$ is nonincreasing, since $D(u_{k+1},u_k)\geq 0$.

(ii) Consider first the case $m=0$. Let $z$ verify $Az=y$ and denote 
$$a=D(z,u_{k+1})+D(u_{k+1},u_k)-D(z,u_k), \,\,\,\,\,\,\xi_k=1+\ln u_k\in \partial f(u_k).$$ 
By using \eqref{optimalitycondition},   one has for all $k\in\mathbb{N}$:
\begin{eqnarray*}
a &=&  cd(z,u_{k+1})-\frac{1}{2}\|Au_{k+1}-y\|^2+cd(u_{k+1},u_k)-\frac{1}{2}\|Au_{k+1}-Au_k\|^2\\
 && - cd(z,u_{k})+\frac{1}{2}\|Au_{k}-y\|^2\\
&=& \langle A^*(Au_k-y), z-u_{k+1}\rangle+\langle Au_{k+1}-y,Au_k-Au_{k+1}\rangle\\
&=& \langle Au_k-y, y-Au_{k+1}\rangle+\langle Au_{k+1}-y,Au_k-Au_{k+1}\rangle\\
&=&-\|y-Au_{k+1}\|^2\leq 0.
\end{eqnarray*}
This implies the typical inequality for a proximal-like method:
\begin{equation}\label{ineq}
\frac{1}{2}\|Au_{k+1}-y\|^2+D(z,u_{k+1})+D(u_{k+1},u_k)\leq D(z,u_k)
\end{equation}
which yields the conclusion.

(iii) Let $m=0$. Inequality \eqref{ineq} leads to
\begin{equation}\label{bound}
\sum_{j=0}^{k}\frac{1}{2}\|Au_{j+1}-y\|^2+D(z,u_{k+1})+\sum_{j=0}^{k}D(u_{j+1},u_j)\leq D(z,u_0),
\end{equation}
which yields 
\begin{equation}\label{mon}
\frac{k+1}{2}\|Au_{k+1}-y\|^2\leq \sum_{j=0}^{k}\frac{1}{2}\|Au_{j+1}-y\|^2,
\end{equation}
since the sequence $\{\|Au_k-y\|^2\}_{k\in\N}$  is monotone.
We show now that $\{f(u_k)\}_{k\in\N}$ is bounded.
To this end, 
due to nonnegativity of $D(z,u_k)$, $k\in\N$, and to  \eqref{ve}, one has
\begin{eqnarray*}
cf(u_k)+\frac{1}{2}\|Au_{k} -y\|^2&\leq & cf(z)+c\langle \xi_k, u_k-z\rangle = cf(z)+\langle w_0+\sum_{j=0}^{k-1}(y-Au_j), Au_k-y\rangle\\
&\leq& cf(z)+ \|w_0\| \|Au_k-y\|+ \sum_{j=0}^{k-1}\frac{1}{2}\|y-Au_j\|^2+\frac{k}{2} \|Au_k-y\|^2.
\end{eqnarray*}
The right hand side is bounded by \eqref{bound} and \eqref{mon}, thus ensuring boundedness of $\{f(u_k)\}_{k\in\N}$. Consequently, there exists a subsequence $\{u_l\}_{l\in\mathbb{N}}$ in $\mbox{dom}\,\partial f$ which is $L^1$-weakly convergent  to some $u\in {\mbox{dom}\, f}$, cf. Lemma \ref{kl_propr}. Then one has $Au_l\rightarrow Au$ weakly in $Y$ and moreover $Au_l\rightarrow y$ in the $Y$-norm since 
\[
\frac{1}{2}\|Au_{k+1}-y\|^2\leq D(z,u_k)-D(z,u_{k+1})\to 0
\]
 due to inequality \eqref{ineq} and to monotonicity of $\{D(z,u_k)\}_{k\in\N}$.  Hence, $u$ satisfies  $Au=y$. 

The proof of the statements above for the case $m=1$  is similar, the main difference being the optimality condition \eqref{optimalitycondition} with $c_k^m$ satisfying $c(\xi_k-\xi_{k+1})=A^*(Au_k-y)+c\ln c_k^m$. In more detail, the term $\langle c\ln c_k^m,z-u_{k+1}\rangle$ vanishes when evaluating $a$ and does not influence further calculations, while other terms containing $c_k^m$ behave similarly  in the remaining argumentation.
\endproof

\end{pro}
Let us consider now the iterative method based on the noisy data, that is
\begin{equation}\label{method0b_noisy}
u_{k+1}\in \mbox{arg}\min_u\left\{\langle Au-\yd,Au_{k}-\yd\rangle +cd(u,u_{k}) +\chi_m(u)\right\}.
\end{equation}

We propose first a discrepancy principle for stopping the algorithm in this case. Before detailing how it works,  denote
\begin{equation}\label{v_noisy}
v_k=w_0+\sum_{j=0}^{k-1}(\yd-Au_j)
\end{equation} 
for $k\geq 1$. Then the optimality condition for \eqref{method0b_noisy}  yields
\begin{equation}\label{ve_noisy}
\xi_k=\xi_0+\sum_{j=0}^{k-1}\ln c_{j}^m+\lambda A^*(\sum_{j=0}^{k-1}(\yd-Au_j))=\sum_{j=0}^{k-1}\ln c_{j}^m+\lambda A^*v_k.
\end{equation}

\begin{pro}\label{conv_noisy} Assume that $A:L^1(\Omega)\rightarrow Y$ is a bounded linear operator which satisfies \eqref{condition}  and such that the operator equation $Au=y$ has a positive solution $z$ verifying $\chi_m(z) = 0$ if $m=1$. Let $\yd\in Y$ be noisy data satisfying $\|y-\yd\|\leq \delta$, for some noise level $\delta$. Let $u_0\in \mbox{dom}\,\partial f$ be an arbitrary starting element with the properties $1+\ln u_0\in \mathcal{R} (A^*)$ and  $\chi_m(u_0) = 0$ if $m=1$. 
 Then 
\begin{enumerate}[(i)]
\item  The residual $\|Au_k-\yd\|$ decreasesmonotonically and the following inequalities hold
\begin{equation}\label{ineq11}
\frac{1}{2}\|\yd-Au_{k+1}\|^2+D(z,u_{k+1})+D(u_{k+1},u_k)\leq \frac{\delta^2}{2}+D(z,u_k),\,\,\,k\in \mathbb{N},
\end{equation}
\begin{equation}\label{ineq12}
\|\yd-Au_{k}\|^2\leq \delta^2+\frac{2D(z,u_0)}{k},\,\,\, k\geq 1.
\end{equation}
\noindent 
\item  The term $D(z,u_k)$ decreases as long as $\|\yd-Au_{k}\|^2>\delta^2.$

\item  The index $k_*(\delta)$  defined by
\begin{equation}\label{aposteriori}
k_*(\delta)=\min\{k\in\mathbb{N}: \|Au_k-y^\delta\|<\sqrt{\tau}\delta\},\,\,\,\tau>1.
\end{equation}
is finite.

\item   There exists a weakly convergent subsequence of $\{u_{k_*(\delta)}\}_\delta$ in $L^1(\Omega)$.
If  $\{{k_*(\delta})\}_\delta$ is unbounded, then each limit point is a solution of $Au=f$.
\end{enumerate}

\proof
We consider only the case $m=0$, since for $m=1$ one can use similar arguments, as explained in the previous proof.
 
\noindent First part of (i) follows by the definition of the iterative procedure.
For proving the remaining inequalities in (i), and (ii), we consider as in the previous proof
\begin{eqnarray*}
a &=&  cd(z,u_{k+1})-\frac{1}{2}\|Au_{k+1}-y\|^2+cd(u_{k+1},u_k)-\frac{1}{2}\|Au_{k+1}-Au_k\|^2\\
& & - cd(z,u_{k})+\frac{1}{2}\|Au_{k}-y\|^2\\
&=& c\langle \xi_k-\xi_{k+1},z-u_{k+1}\rangle+\langle Au_{k+1}-y,Au_k-Au_{k+1}\rangle\\
&=& \langle Au_k-\yd, y-Au_{k+1}\rangle+\langle Au_{k+1}-y,Au_k-Au_{k+1}\rangle\\
&=&  \langle y- Au_{k+1}, Au_{k+1}-\yd\rangle =-\|\yd-Au_{k+1}\|^2+ \langle y- \yd, Au_{k+1}-\yd\rangle\\
&\leq& -\|\yd-Au_{k+1}\|^2+\frac{\delta^2}{2}+\frac{1}{2}\|\yd-Au_{k+1}\|^2\\
&=& -\frac{1}{2}\|\yd-Au_{k+1}\|^2+\frac{\delta^2}{2}.
\end{eqnarray*}
Inequality \eqref{ineq12} can be obtained by writing \eqref{ineq11} for $k=0,...,n-1$ und calculating the telescope sum:
\begin{equation*}
 \frac{n \|\yd-Au_{n}\|^2}{2}\leq\frac{1}{2} \sum_{k=0}^{n-1} \|\yd-Au_{k+1}\|^2+D(z,u_n)\leq D(z,u_0)+\frac{\delta^2n}{2}.
\end{equation*}
Moreover, (ii) follows from \eqref{ineq11} by neglecting $D(u_{k+1},u_k)$.

\noindent(iii) follows from \eqref{ineq12} and the definition of $k_*(\delta)$:
\begin{equation}\label{bound_noisy}
\frac{k_*(\delta)\tau\delta^2}{2}\leq \frac{1}{2} \sum_{k=0}^{k_*(\delta)-1} \|\yd-Au_{k}\|^2+D(z,u_{k_*(\delta)})\leq D(z,u_0)+\frac{\delta^2k_*(\delta)}{2},
\end{equation}
which implies 
\begin{equation}\label{index}
k_*(\delta)\leq\frac{2D(z,u_0)}{(\tau-1)\delta^2}.
\end{equation}

\noindent (iv)  can be shown similarly to Proposition \ref{conv_exact} (iii). Due to nonnegativity of $D(z,u_k)$ for any $k\in\N$ and to \eqref{ve_noisy}, one has
\begin{eqnarray*}
cf(u_k)+\frac{1}{2}\|Au_{k} -y\|^2&\leq & cf(z)+c\langle \xi_k, u_k-z\rangle\\
&=& cf(z)+\langle w_0+\sum_{j=0}^{k-1}(\yd-Au_j), Au_k-y\rangle\\
&\leq& cf(z)+ \|w_0\| \|Au_k-y\|+ \sum_{j=0}^{k-1}\|\yd-Au_j\|\| Au_k-y\|\\
&\leq& cf(z)+ \|w_0\| \|Au_k-\yd\|+ \delta \|w_0\|+\delta \sum_{j=0}^{k-1}\|\yd-Au_j\|\\
&+&\sum_{j=0}^{k-1}\|\yd-Au_j\|\| Au_k-\yd\|\\
&\leq& cf(z)+ \|w_0\| \|Au_k-\yd\|+ \delta \|w_0\| +\sum_{j=0}^{k-1}\frac{1}{2}\|y-Au_j\|^2+ \frac{k\delta^2}{2}\\ 
&+& \sum_{j=0}^{k-1}\frac{1}{2}\|y-Au_j\|^2+\frac{k}{2} \|Au_k-\yd\|^2.
\end{eqnarray*}
The right hand side written for $k=k_*(\delta)$ is bounded by  \eqref{bound_noisy}, \eqref{aposteriori}, the monotonicity of the residual and by \eqref{index}, thus ensuring boundedness of $\{f(u_{k_*(\delta)}\}_{\delta>0}$ for $\delta$ small enough. The conclusion follows then as in the proof of Proposition \ref{conv_exact} (iii). \endproof

\end{pro}

A convergence result can be established also in case of an a priori stopping rule with $k_*(\delta)\sim\frac{1}{\delta}$ by following the lines of Proposition \ref{conv_noisy} (iv).

\begin{pro} Assume that $A:L^1(\Omega)\rightarrow Y$ is a bounded linear operator which satisfies \eqref{condition}  and such that the operator equation $Au=y$ has a positive solution $z$ verifying $\chi_m(z) = 0$ if $m=1$. Let $\yd\in Y$ be noisy data satisfying $\|y-\yd\|\leq \delta$, for some noise level $\delta$. Let $u_0\in \mbox{dom}\,\partial f$ be an arbitrary starting element with the properties $1+\ln u_0\in \mathcal{R} (A^*)$ and  $\chi_m(u_0) = 0$ if $m=1$.  Let the stopping index $k_*(\delta)$ be chosen of order
$1/{\delta}$.   Then $\{f(u_{k_*(\delta)})\}_{\delta}$ is  bounded and
hence, as $\delta \rightarrow 0$, there exists a weakly convergent
subsequence $\{u_{k(\delta_n)}\}_n$ in $L^1(\Omega)$ whose limit is a solution of
$Au =y$. Moreover, if the solution of the equation is unique, then
$\{u_{k_*(\delta)}\}_{\delta>0}$ converges weakly to the solution as
${\delta\rightarrow 0}$.

\end{pro}

\section{Error estimates}

In this section we derive error estimates under a specific source condition (on a solution) for the entropy type penalty. We proceed first with the case of exact data on the right-hand side of the operator equation and then with the noisy data case, by employing an a priori rule for stopping the algorithm.

\subsection{Exact data case}

\begin{pro}\label{pro:exact_data} Assume that $A:L^1(\Omega)\rightarrow Y$ is a bounded linear operator which satisfies \eqref{condition}  and such that the operator equation $Au=y$ has a positive solution $z$ verifying $\chi_m(z) = 0$ if $m=1$. Let $u_0\in \mbox{dom}\,\partial f$ be an arbitrary starting element with the properties $1+\ln u_0\in \mathcal{R} (A^*)$ and  $\chi_m(u_0) = 0$ if $m=1$. 
Additionally, let the following source condition hold: 
\begin{equation}\label{sc}
1+\ln z\in \mathcal{R}(A^*).
\end{equation}
Then one has 
\begin{equation}\label{estimate}
d(z,u_k)=O(1/k).
\end{equation}
Moreover, $\|u_k-z\|_1=O(1/\sqrt{k})$ if $m=1$.

\proof We consider only the case $m=0$ (similar arguments for the other case).

First, we symmetrize  $D$ by considering $D^s(x,y)=D(x,y)+D(y,x)$. Let $\xi=1+\ln z=\lambda A^*v$ for some $v\in Y$.

One can use similar techniques as in \cite{BurResHe} for deriving the announced error estimates, by carefully dealing with  the setting of the  $D$ distance penalty. Based on \eqref{v}, one has
\begin{eqnarray*}
D^s(u_k,z)&=&c\langle \xi_k-\xi,u_k-z\rangle-\|Au_k-Az\|^2\\
&=&\langle A^*v_k-A^*v,u_k-z\rangle-\|Au_k-y\|^2\\
&=&\langle v_k-v,Au_k-y\rangle-\|Au_k-y\|^2\\
&=&\langle v_k-v,v_k-v_{k+1}\rangle-\|Au_k-y\|^2\\
&=&\frac{1}{2}\| v_k-v\|^2-\frac{1}{2}\| v_{k+1}-v\|^2+\frac{1}{2}\| v_{k+1}-v_k\|^2-\|Au_k-y\|^2\\
&=&\frac{1}{2}\| v_k-v\|^2-\frac{1}{2}\| v_{k+1}-v\|^2-\frac{1}{2}\|Au_k-y\|^2.
\end{eqnarray*}
By writing the last inequality  also  for $k-1, k-2,...,1$, by summing up and by combining  with monotonicity of  $\{D(z,u_k)\}$, one obtains
$$kD(z,u_k)\leq \sum_{j=1}^{k}D(z,u_j)\leq \sum_{j=1}^{k}D^s(u_j,z)\leq \frac{1}{2}\| v_1-v\|^2-\frac{1}{2}\| v_{k+1}-v\|^2-\frac{1}{2}\sum_{j=1}^k\|Au_j-y\|^2$$
and thus, due to \eqref{mon},
$$d(z,u_k)\leq \frac{\lambda}{2k}\| v_1-v\|^2$$
holds. The announced convergence rate in the $L^1$-norm holds in case $m=1$ by Lemma \ref{entropy_propr} (v).
\endproof

\end{pro}

\subsection{Noisy data case}

\begin{pro} Assume that $A:L^1(\Omega)\rightarrow Y$ is a bounded linear operator which satisfies \eqref{condition}  and such that the operator equation $Au=y$ has a positive solution $z$ verifying $\chi_m(z) = 0$ if $m=1$. Let $u_0\in \mbox{dom}\,\partial f$ be an arbitrary starting element with the properties $1+\ln u_0\in \mathcal{R} (A^*)$ and  $\chi_m(u_0) = 0$ if $m=1$. Let $\yd\in Y$ be noisy data satisfying $\|y-\yd\|\leq \delta$, for some noise level $\delta$. 
Let the stopping index $k_*(\delta)$ be chosen of order $1/\delta$.
 and let the  source condition \eqref{sc} hold.
Then one has 
\begin{equation}\label{estimate_noisy}
d(z,u_{k_*(\delta)})=O(\delta).
\end{equation}
Moreover, $\|u_{k_*(\delta)}-z\|_1=O(\sqrt{\delta})$ if $m=1$.

\proof
Note that  $c\xi-A^*Az=A^*q$  for  $q=v-A^*z$. With this notation,
one can show the following estimate as in  Theorem 4.3 in \cite{BurResHe}:
$$D(z,u_k)\leq \frac{\|q\|^2}{2k}+\delta\|q\|+\delta^2k, \,\,\,\forall k\in\N.$$
Then one has
\begin{eqnarray*}
cd(z,u_k)\leq \frac{\|q\|^2}{2k}+\delta\|q\|+\delta^2k+\frac{\|Au_k-y\|^2}{2}\leq  \frac{\|q\|^2}{2k}+\delta\|q\|+\delta^2k+\|Au_k-\yd\|^2+\delta^2
\end{eqnarray*}
which yields \eqref{estimate_noisy} when written for $k=k_*(\delta)$, due to \eqref{ineq12}.
\endproof\smallskip

\noindent Establishing convergence rates by means of a discrepancy rule remains an open issue.
\end{pro}

\section{General data fidelities}
Before we conclude with numerical examples, we want to emphasize that Problem \eqref{method0} can easily be generalised to
\begin{align}
u_{k + 1} \in \arg\min_u \left\{ F_{y^\delta}(Au) + c d(u, u_k) -  g(u, u_k)) + \chi_m(u) \right\} \, .\label{eq:methodgen}
\end{align}
Here $F_{y^\delta}:Y  \rightarrow [0, +\infty)$ is a more general data fidelity term that is assumed to be convex and Fr\'{e}chet-differentiable and $g: L^1(\Omega) \times L^1(\Omega) \rightarrow [0,+\infty)$ is the Bregman distance with respect to the function $F_{y^\delta}$, i.e.
\begin{align*}
g(u, v) := F_{y^\delta}(Au) - F_{y^\delta}(Av) - \langle A^\ast F_{y^\delta}^\prime(Av), u - v \rangle \, .
\end{align*}
Note that \eqref{eq:methodgen} is an instance of the Bregman proximal method \cite{CensorZenios92,Gutman}. The update for \eqref{eq:methodgen} can be written, in analogy to \eqref{entropic},  as
\begin{align}
u_{k + 1}^m = u_k^m c_k^m e^{-\lambda A^\ast F_{y^\delta}^\prime(Au_k)} \, ,\label{eq:genentropic}
\end{align}
for $\lambda = 1/c$ and $m \in \{0, 1\}$. We want to emphasise that a more general data fidelity term that satisfies the assumptions mentioned above together with 
\begin{align}
c d(u, v) -  g(u, v) \geq 0 \, ,\label{eq:gensmoothness}
\end{align}
for all $u, v \in dom \, f$, is no restriction in terms of Fej\'{e}r-monotonicity. In analogy to \cite[Lemma 6.11]{BenningBurger18} we can conclude
\begin{align*}
D(z, u_{k + 1}) \leq D(z, u_k)
\end{align*}
for all $k < k_\ast(\delta)$, with $k_\ast(\delta)$ chosen according to a modified version of \eqref{aposteriori} that reads as
\begin{align}
k_*(\delta)=\min\{k\in\mathbb{N}: F_{y^\delta}(Au_k) < \delta \} \, .\label{eq:modaposteriori}
\end{align}
However, we can also derive a monotonicity result for $d$ directly. First of all we observe that \eqref{eq:gensmoothness} implies
\begin{align*}
-\langle A^\ast F_{y^\delta}^\prime(Au_k), u_{k + 1} - u_k \rangle - c d(u_{k + 1}, u_k) \leq F_{y^\delta}(Au_k) - F_{y^\delta}(Au_{k + 1}) \, .
\end{align*}
Inserting \eqref{eq:genentropic} into the inequality above then yields
\begin{align*}
c \langle \ln u_{k + 1} - \ln u_k, u_{k + 1} - u_k \rangle - c d(u_{k + 1}, u_k) \leq F_{y^\delta}(Au_k) - F_{y^\delta}(Au_{k + 1}) + c \langle \ln c_k^m, u_{k + 1} - u_k \rangle \, .
\end{align*}
As mentioned earlier in Section \ref{sec:convergence}, we either have $\ln(c_k^m) = 0$ for $m = 0$, or orthogonality of $\ln(c_k^m)$ to all functions of the form $w - v$ with $\chi_m(w) = \chi_m(v) = 0$ for $m = 1$. Since $d(u_{k + 1}, u_k) + d(u_k, u_{k + 1}) = \langle \ln u_{k + 1} - \ln u_k, u_{k + 1} - u_k \rangle$, we therefore estimate
\begin{align*}
c d(u_k, u_{k + 1}) \leq F_{y^\delta}(Au_k) - F_{y^\delta}(Au_{k + 1}) \, .
\end{align*} 
With the three-point identity we then observe
\begin{align*}
d(z, u_{k + 1}) - d(z, u_k) {} = {} &- \langle \ln u_{k + 1} - \ln u_k, z - u_{k + 1} \rangle - d(u_{k + 1}, u_k) \\
{} = {} &- \langle \ln u_{k + 1} - \ln u_k, z - u_k \rangle + \langle \ln u_{k + 1} - \ln u_k, u_{k + 1} - u_k \rangle \\
&- d(u_{k + 1}, u_k) \\
{} = {} &- \langle \ln u_{k + 1} - \ln(u_k), z - u_k \rangle + d(u_k, u_{k + 1}) \\
{} = {} &\langle \lambda A^\ast F_{y^\delta}^\prime(Au_k), z - u_k \rangle + d(u_k, u_{k + 1}) \\
{} \leq {} &\lambda \left( F_{y^\delta}(Az) - F_{y^\delta}(Au_k) \right) + d(u_k, u_{k + 1}) \\
{} \leq {} & \lambda \left( F_{y^\delta}(Az) - F_{y^\delta}(Au_{k + 1}) \right) \, .
\end{align*}
Together with \eqref{eq:modaposteriori} we can then conclude
\begin{align*}
d(z, u_{k + 1}) < d(z, u_k)
\end{align*}
for $k < k^\ast(\delta)$.


\section{Example Problems}

We finally discuss several types of problems that satisfy the conditions used in the analysis and present numerical illustrations for some of these situations. 

\subsection{Integral Equations}

Let $\Omega \subset \R^d$ and $\tilde \Omega \subset \R^{\tilde d}$ be open and bounded sets and let $k \in L^\infty(\tilde \Omega \times \Omega)$.
Then the integral operator
\begin{equation}
	A: L^1(\Omega) \rightarrow L^2(\tilde \Omega) , \qquad u \mapsto \int_\Omega k(\cdot,y) u(y) ~dy
\end{equation}
is a well-defined and bounded linear operator. Thus, the convergence analysis is applicable due to Lemma  \ref{continuity}. 

We mention that in the case of $k$ being a nonnegative function, and hence $A$ and $A^*$ preserving nonnegativity, standard schemes preserving nonnegativity are available. In particular for $k$ including negative entries, the entropic Landweber scheme offers a straightforward alternative, since it does not depend on the positivity preservation of $A$ respectively its adjoint. For comparison we consider the EM-Algorithm
$$ u_{k+1} = \frac{u_k}{A^* 1} A^* \left( \frac{y}{Au_k} \right)$$ and 
the projected Landweber iteration
$$ u_{k+1} = \left( u_k - \tau A^*(Au_k -y) \right)_+. $$

\noindent We implement the forward operator by discretization of $u$ on a uniform grid and a trapezoidal rule for integration. We use the following examples of kernels and initial values, all on $\Omega = (0,1)$, the first two being standard test examples used in the literature on maximum entropy methods (cf. \cite{amato})
\begin{enumerate}

\item Kernel $k_1(x,y)=e^{xy}$, exact solution $z_1(x)=e^{-\frac{x^2}{2\sigma^2}}$

\item Kernel $k_2(x,y)=3e^{-\frac{(x-y)^2}{0.04}}$, exact solution 
$$z_2(x)=1-0.9 e^{-\frac{(x-0.1)^2}{2\sigma^2}}-0.3 e^{-\frac{(x-0.3)^2}{2\sigma^2}}
-0.5 e^{-\frac{(x-0.5)^2}{2\sigma^2}}-0.2 e^{-\frac{(x-0.7)^2}{2\sigma^2}}-0.7 e^{-\frac{(x-0.9)^2}{2\sigma^2}}.$$  

\item Kernel $k_3(x,y)=1$ if $x\geq y$ and $k_3(x,y) = 0$ else, exact solution $z_3(x)=e^{-\frac{x^2}{2\sigma^2}}$ .
\end{enumerate}

\begin{figure} 
	\centering
	\includegraphics[width=.48\textwidth]{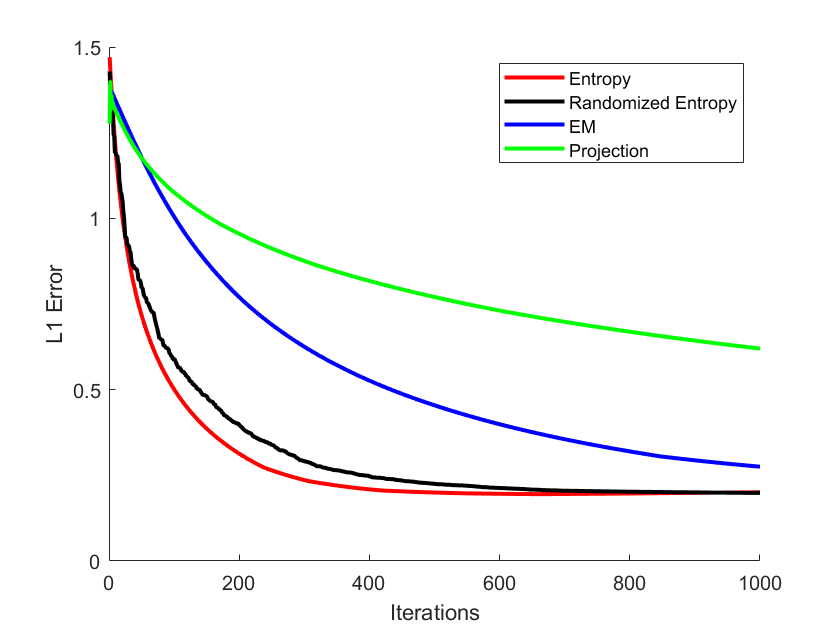} 
	\includegraphics[width=.48\textwidth]{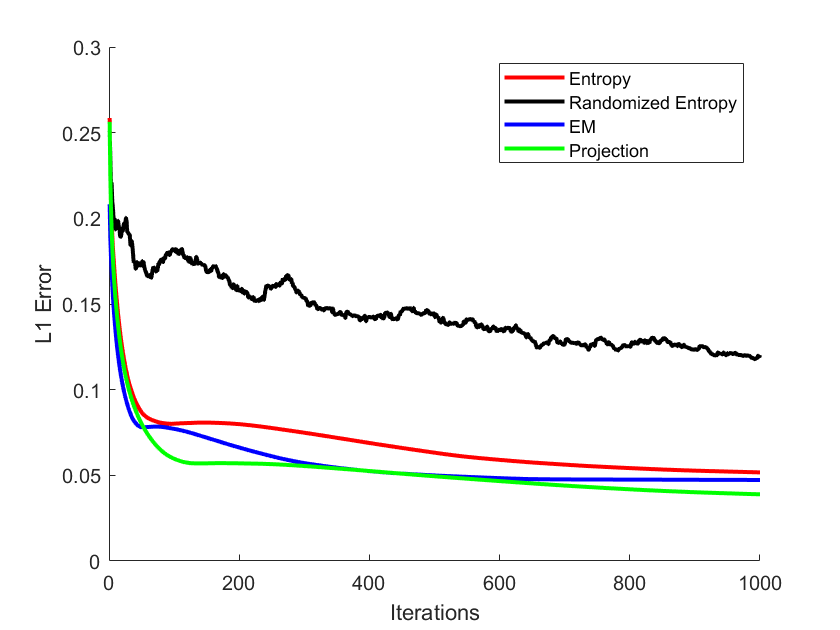} 
	
	\includegraphics[width=.48\textwidth]{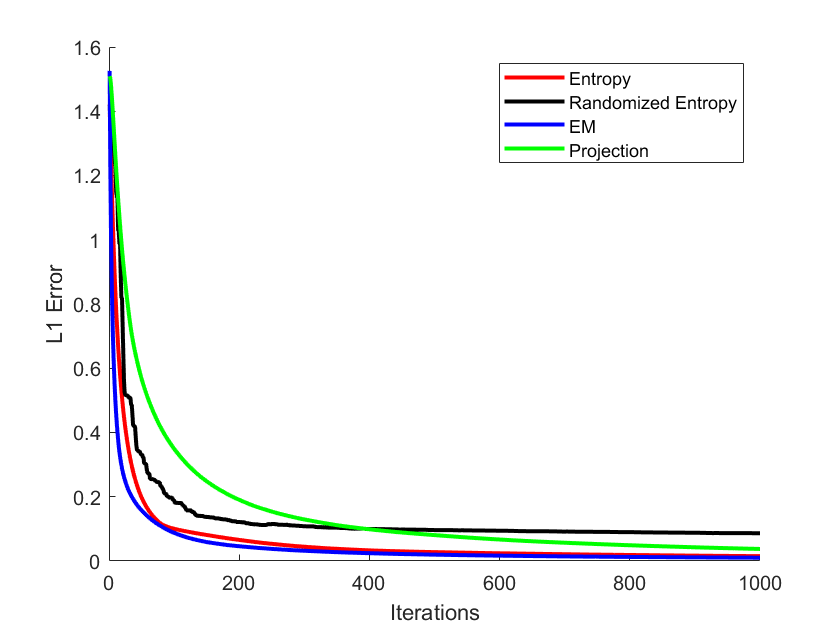} 
	\caption{Results of the three test cases,  $L^1$ error plotted vs. iteration number.  \label{resultsfigure} }
\end{figure}

In all examples, we have chosen  $\sigma^2 =0.01$ and a constant intial value $u^0$. In order to illustrate the behaviour of the iteration methods we plot the the error $\Vert u^k - z \Vert_{L^1}$ vs. the iteration number $k$ in Figure \ref{resultsfigure}. 

We observe that the entropic projection is at least competitive to the other schemes in all examples, it outperforms the EM and projection method in the first example, which is a combination of severe ill-posedness with an exact solution having many entries close to zero (which is a particularly difficult case for the EM algorithm). 

In the second case, again severely ill-posed, the projected Landweber iteration performs better, mainly due to  the strong initial decrease when the solution is positive and no projection is applied. 

The third case corresponding to numerical differentiation, i.e. a very mildly ill-posed problem, is characterized by fast convergence of the schemes, but again the projection method converges significantly slower. For comparison we also include the stochastic version of the entropic projection method, with only one equation used in each iteration step, hence a highly efficient computation. That is,  the operator $A$ is divided in $M$ blocks $A=(A_1,A_2,\ldots,A_M)^T$, and  the data $y$ are partitioned in the same way: $y=(y_1,y_2,\ldots,y_M)$.
With $J(k)$ a discrete uniform random variable in $\{1,\ldots,M\}$, we compute the iterates
$$ u_{k+1}= u_{k} c_{k}^m  e^{\lambda M A_{J(k)}^*(y_{J(k)}-A_{J(k)} u_{k})},\,\,\,k\in\N.$$
 The initial convergence curve is similar to the other method, with much lower computational effort, then  the asymptotic convergence close to the exact solution becomes significantly slower. Hence, it might be very attractive to use the stochastic version at least for the first phase of the reconstruction.

\begin{figure}[t]
\begin{center}
\begin{subfigure}[b]{0.49\textwidth}
\centering
\includegraphics[width=\textwidth]{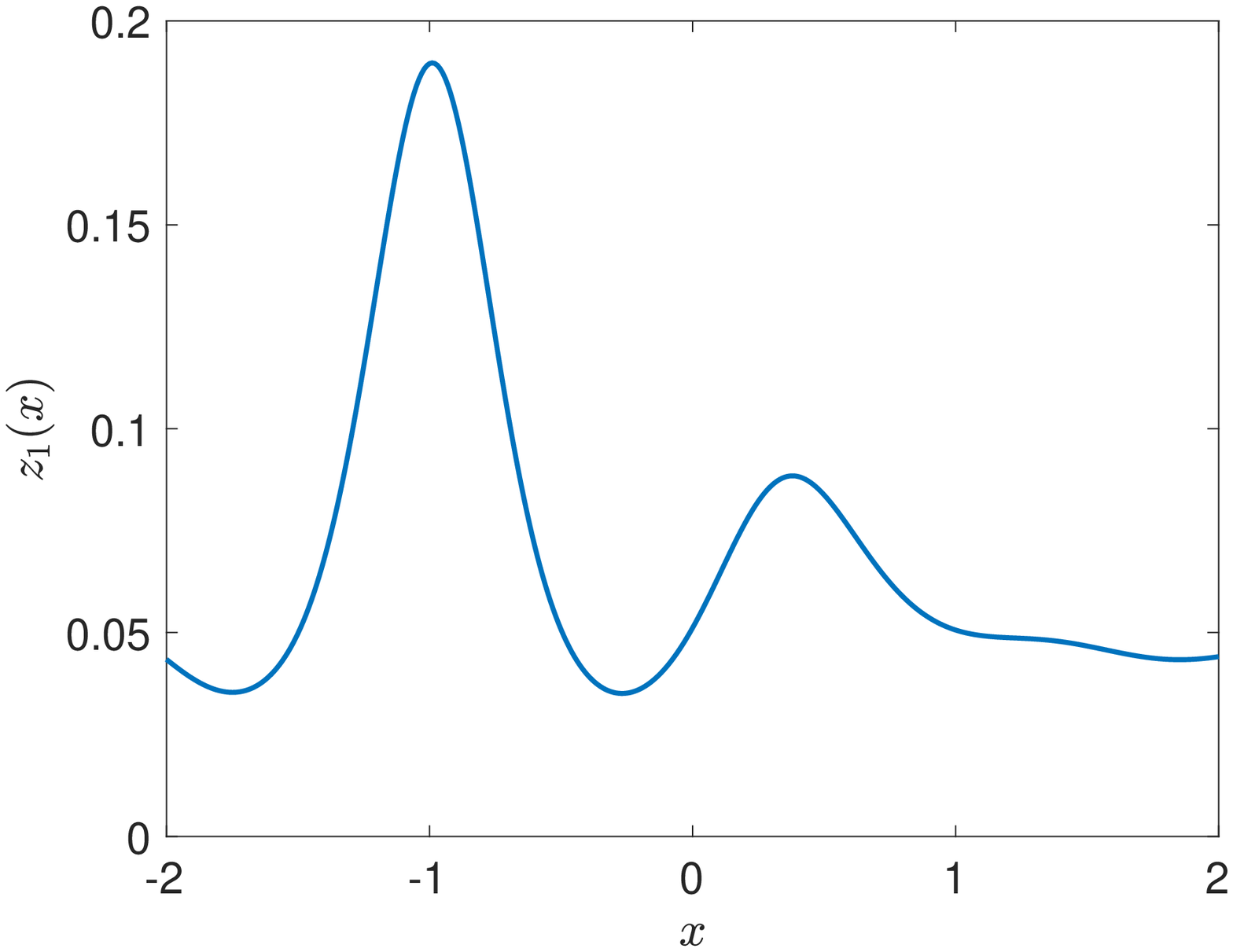}
\caption{$z_1$}
\label{subfig:sampling_gt1}
\end{subfigure}
\begin{subfigure}[b]{0.49\textwidth}
\centering
\includegraphics[width=\textwidth]{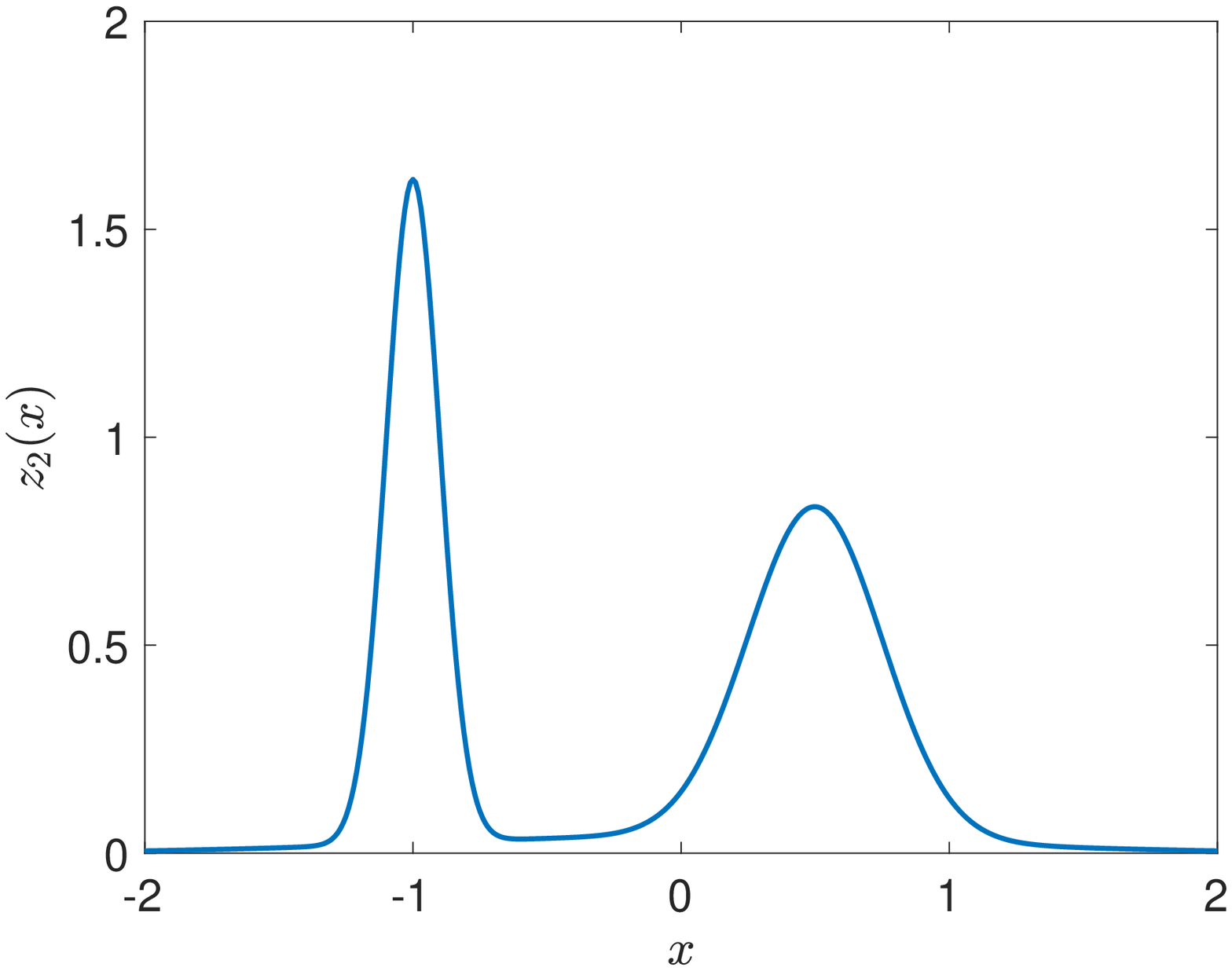}
\caption{$z_2$}
\label{subfig:sampling_gt2}
\end{subfigure}
\end{center}
\caption{The two ground truth functions $z_1$ and $z_2$ as defined in Section \ref{sec:discrete}. Note that $z_1$ by construction satisfies \eqref{sc} while $z_2$ does not satisfy \eqref{sc}.}
\label{fig:sampling_ground_truth}
\end{figure}

\subsection{Discrete sampling of continuous probability densities}\label{sec:discrete}
Suppose that our forward operator is the Fourier integral of a real-valued function evaluated at discrete samples $\xi_1, \ldots, \xi_n$ on a compact domain $\Omega \subset \R^d$, i.e.
\begin{align*}
A : L^1(\Omega) \rightarrow \C^n \, , \qquad u \mapsto \left((2 \pi)^{-\frac{d}{2}} \int_{\Omega} u(x) \, e^{-i x \cdot \xi_j} \, dx \,\right)_{1\leq j\leq n} .
\end{align*}
Then the adjoint operator $A^\ast$ that satisfies $\sum_{j = 1}^n (Au)_j \overline{v_j} = \int_{\Omega} u(x) \overline{(A^\ast v)(x)} dx = \int_{\Omega} u(x) (A^\ast v)(x) dx$ is given as
\begin{align*}
A^\ast : \C^n \rightarrow L^\infty(\Omega) \, , \qquad v \mapsto \text{Re}\left( (2 \pi)^{-\frac{d}{2}} \sum_{j = 1}^n v_j \, e^{i x \cdot \xi_j} \right) \, ,
\end{align*}
where $\text{Re}$ denotes the real part of a complex function. For this choice of $A$ the iterates of \eqref{method0b} read
\begin{align*}
u_{k + 1}(x) = c_k^m u_k(x) e^{\lambda \text{Re}\left( (2 \pi)^{-\frac{d}{2}} \sum_{j = 1}^n \left( y_j - (2 \pi)^{-\frac{d}{2}} \int_{\Omega} u_k(t) e^{-i t \cdot \xi_j} \, dt  \right) \, e^{i x \cdot \xi_j} \right) } \, ,
\end{align*}
where $c_k^m$ is defined as in \eqref{entropic}. Note that this update can also be written as
\begin{align*}
u_k(x) &= c_{k - 1}^m \dots c_0^m u_0(x) e^{\lambda \text{Re}\left( (2 \pi)^{-\frac{d}{2}} \sum_{j = 1}^n \left( k \, y_j - \tilde{y}_{j k}  \right) \, e^{i x \cdot \xi_j} \right) } \, ,
\intertext{for}
\tilde{y}_{j k} &:= (2 \pi)^{-\frac{d}{2}} \sum_{l = 0}^{k - 1} \int_{\Omega} u_l(t) e^{-i t \cdot \xi_j} \, dt \, .
\end{align*}
Due to $\tilde{y}_{j (k + 1)} = \tilde{y}_{j k} + (2 \pi)^{-\frac{d}{2}} \int_{\R^d} u_k(t) e^{-i t \cdot \xi_j} \, dt$, this formulation has the advantage that the numerical costs for evaluating the integrals remains constant. 

\begin{figure}[t]
\begin{subfigure}[b]{0.49\textwidth}
\centering
\includegraphics[width=\textwidth]{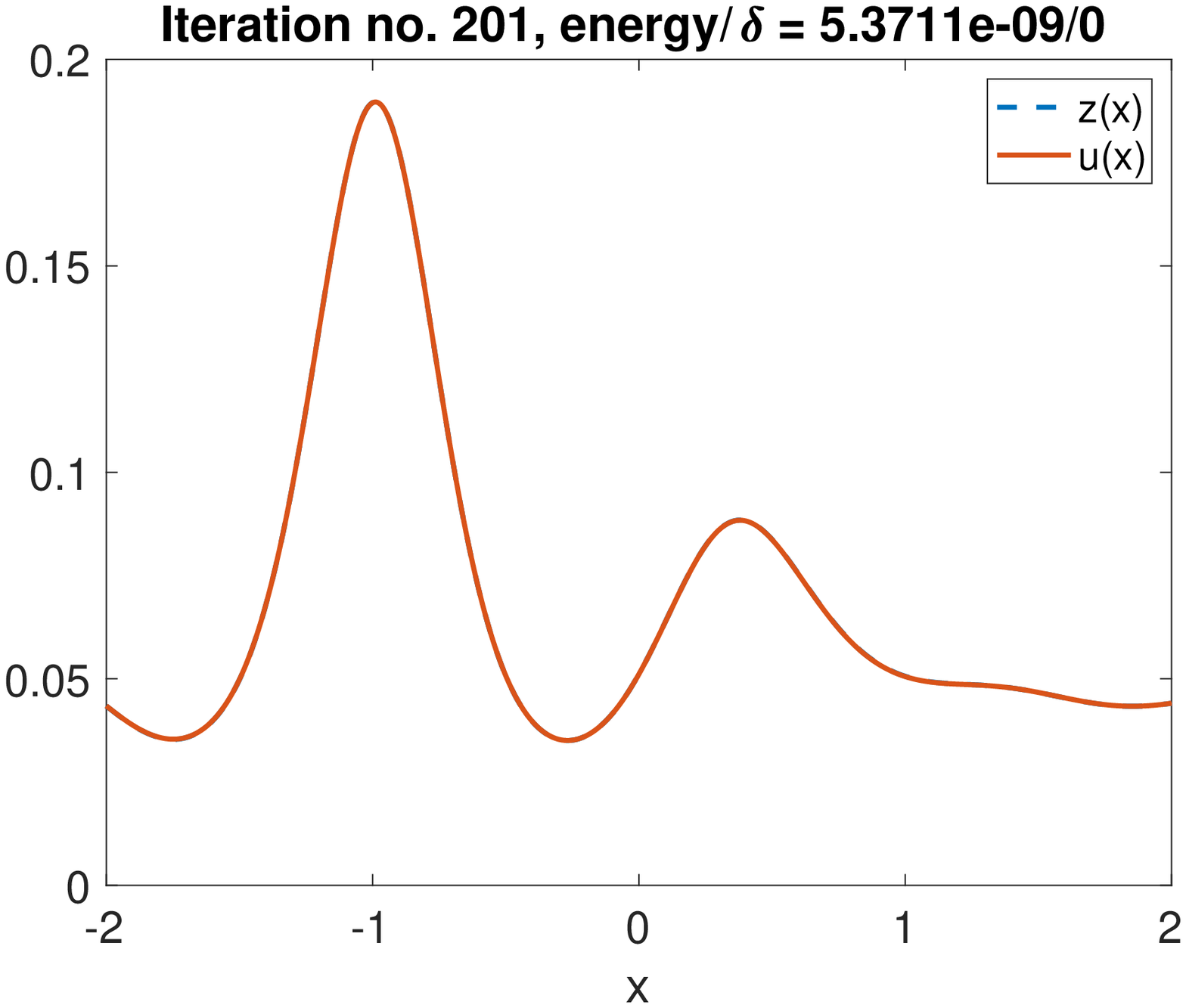}
\caption{$u^{k^\ast}$}
\label{subfig:results_gt1_1}
\end{subfigure}
\begin{subfigure}[b]{0.49\textwidth}
\centering
\includegraphics[width=\textwidth]{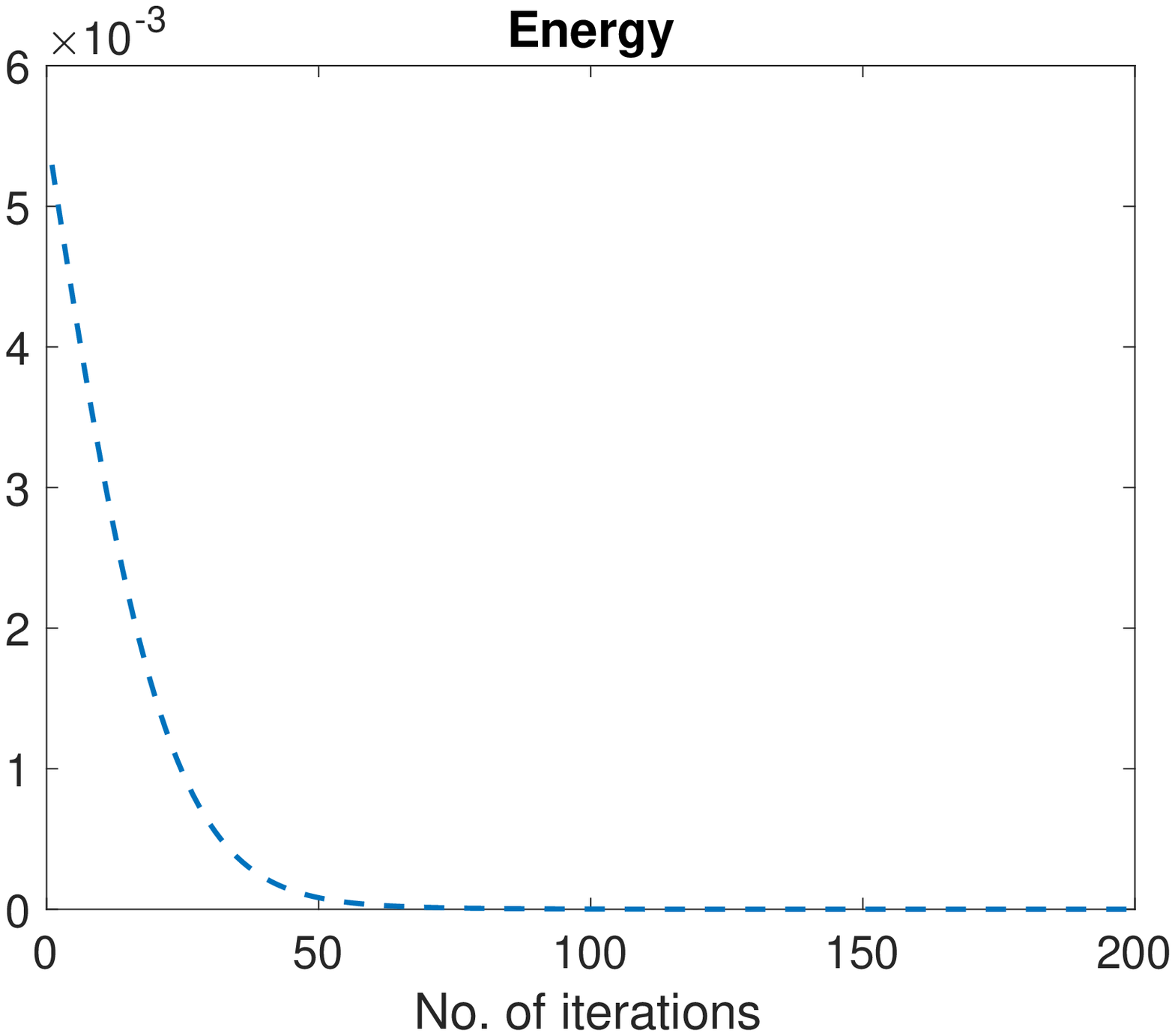}
\caption{Energy decrease of objective}
\label{subfig:results_gt1_2}
\end{subfigure}\\
\begin{subfigure}[b]{0.49\textwidth}
\centering
\includegraphics[width=\textwidth]{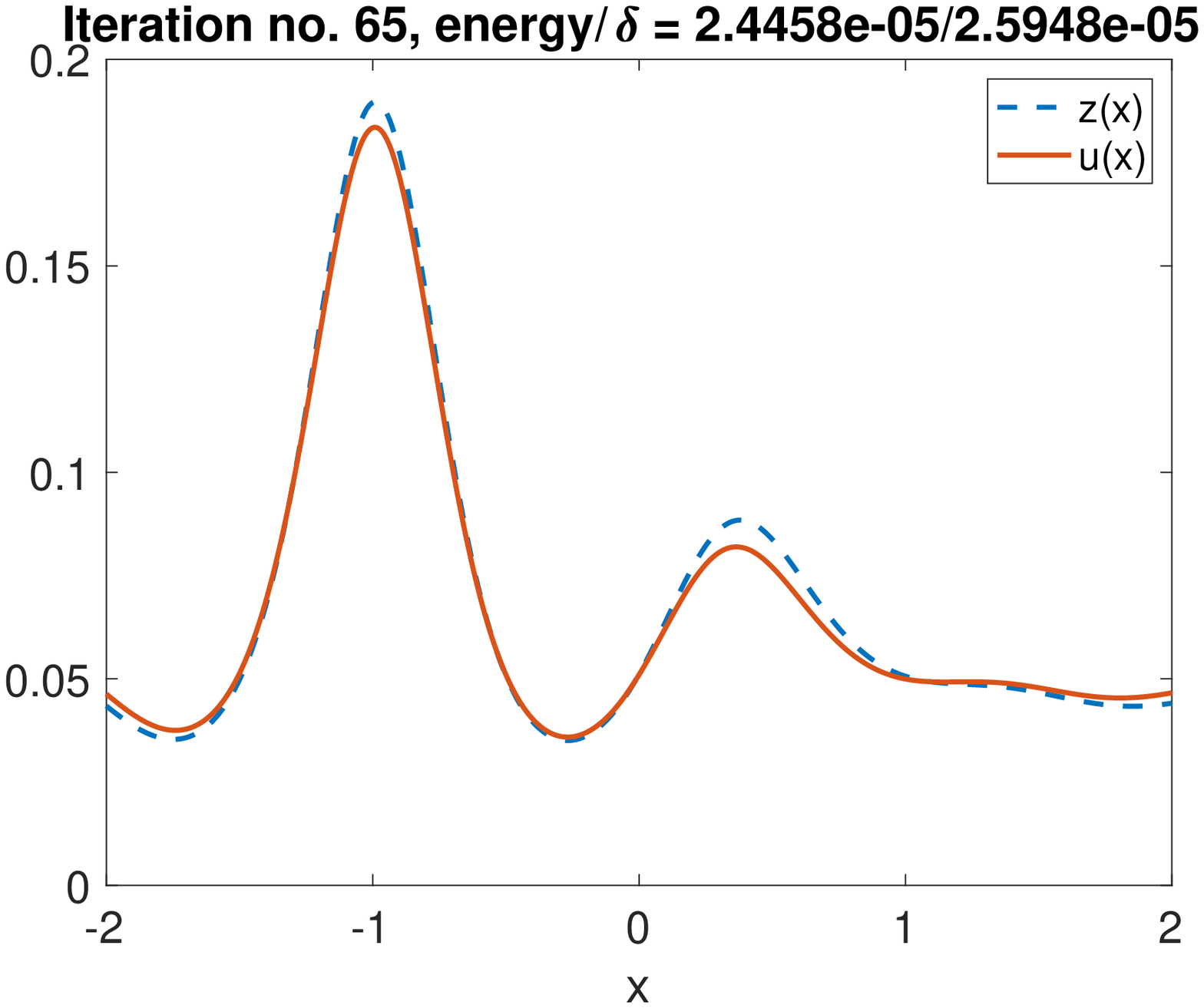}
\caption{$u^{k^\ast}$}
\label{subfig:results_gt1_3}
\end{subfigure}
\begin{subfigure}[b]{0.49\textwidth}
\centering
\includegraphics[width=\textwidth]{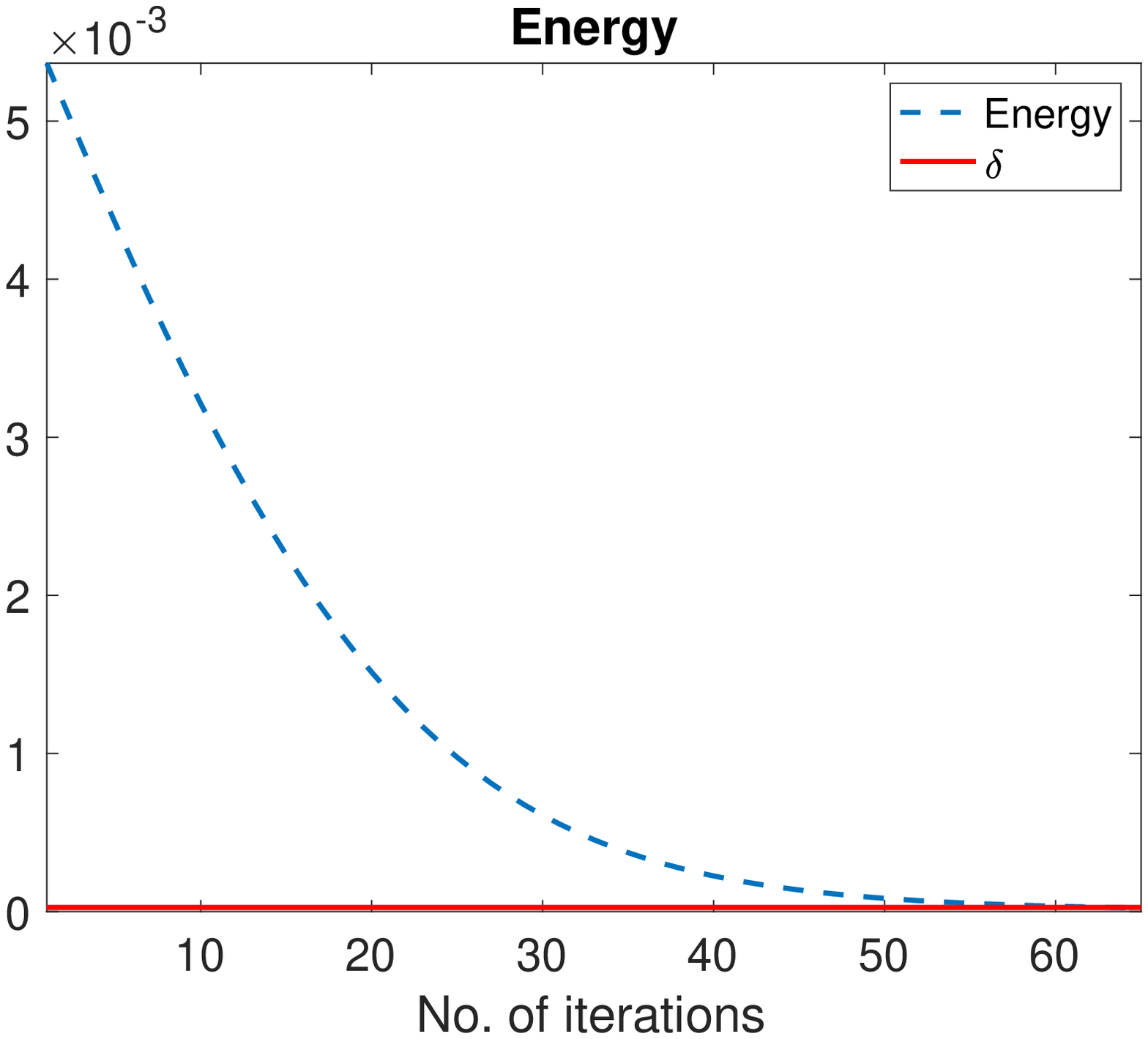}
\caption{Energy decrease of objective}
\label{subfig:results_gt1_4}
\end{subfigure}
\caption{Results of algorithm \eqref{method0b} for data measurements of the form \eqref{eq:measurements} based on the underlying function $z_1$ as visualised in Figure \ref{subfig:sampling_gt1}. Figure \ref{subfig:results_gt1_1} shows the ground truth function $z_1$ and $u^{k^\ast}$ under the assumption of $\sigma = 0$ in \eqref{eq:measurements} , whereas Figure \ref{subfig:results_gt1_2} visualises the monotonic energy decrease over the course of the iterations. Figure \ref{subfig:results_gt1_3} and Figure \ref{subfig:results_gt1_4} show the same results under the assumption of $\sigma = 1/50$ in \eqref{eq:measurements}.}
\label{fig:sampling_results1}
\end{figure}

In the following we consider a one-dimensional setting ($d = 1$) with $\Omega = [-a, a]$ for $a = 10$, where we measure $n = 16$ samples $\{ y_j \}_{j = 1}^n$ of the Fourier integral for coordinates $\xi_j = (2 \pi (j - 1))/n$, $j = 1, \ldots, n$. We assume that these measurements are of the form
\begin{align}
y_j = \frac{1}{\sqrt{2 \pi}} \int_{-10}^{10} z(t) \, e^{- i t \xi_j} \, dt + n_j \, ,\label{eq:measurements}
\end{align}
for a function $z \in L^1([-10, 10])$ and where $n_j \in \mathcal{N}(0, \sigma^2)$ are normal-distributed random variables with mean zero and variance $\sigma^2$, for all $j \in \{1, \ldots, n\}$. We consider numerical experiments for two choices of $z$. The first choice is the following Gau\ss ian-mixture model,
\begin{align*}
\tilde{z}(x) := \sum_{l = 1}^3 c_l \, g(x, \mu_l, \sigma_l) \, ,
\end{align*}
that is constructed as a linear combination of three normalised Gau\ss ians, i.e.
\begin{align*}
g(x, \mu, \sigma) := \frac{1}{\sqrt{2 \pi \sigma^2}} e^{- \frac{(x - \mu)^2}{2 \sigma^2}} \, .
\end{align*}

\begin{figure}[t]
\begin{subfigure}[b]{0.49\textwidth}
\centering
\includegraphics[width=\textwidth]{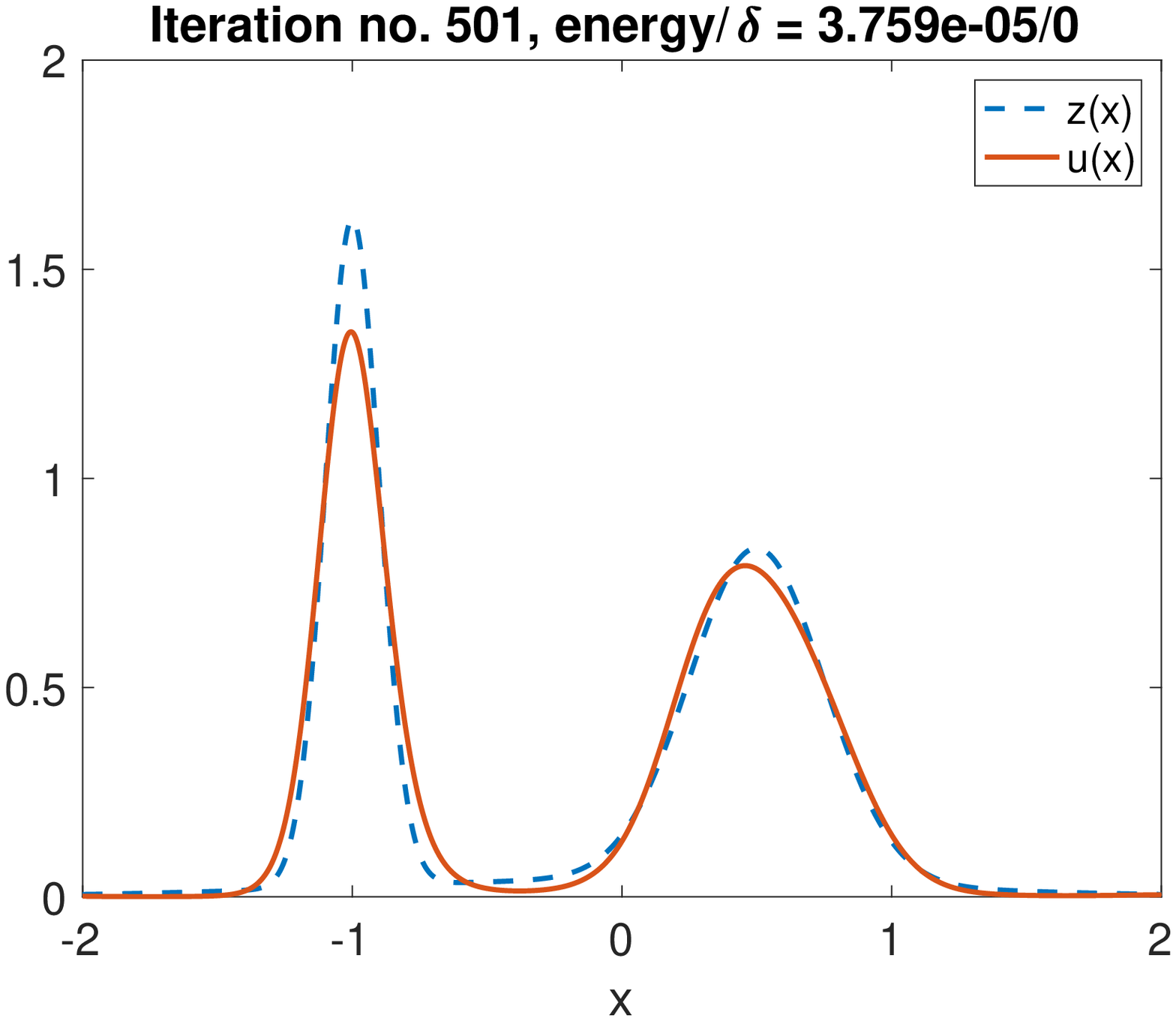}
\caption{$u^{k^\ast}$}
\label{subfig:results_gt2_1}
\end{subfigure}
\begin{subfigure}[b]{0.49\textwidth}
\centering
\includegraphics[width=\textwidth]{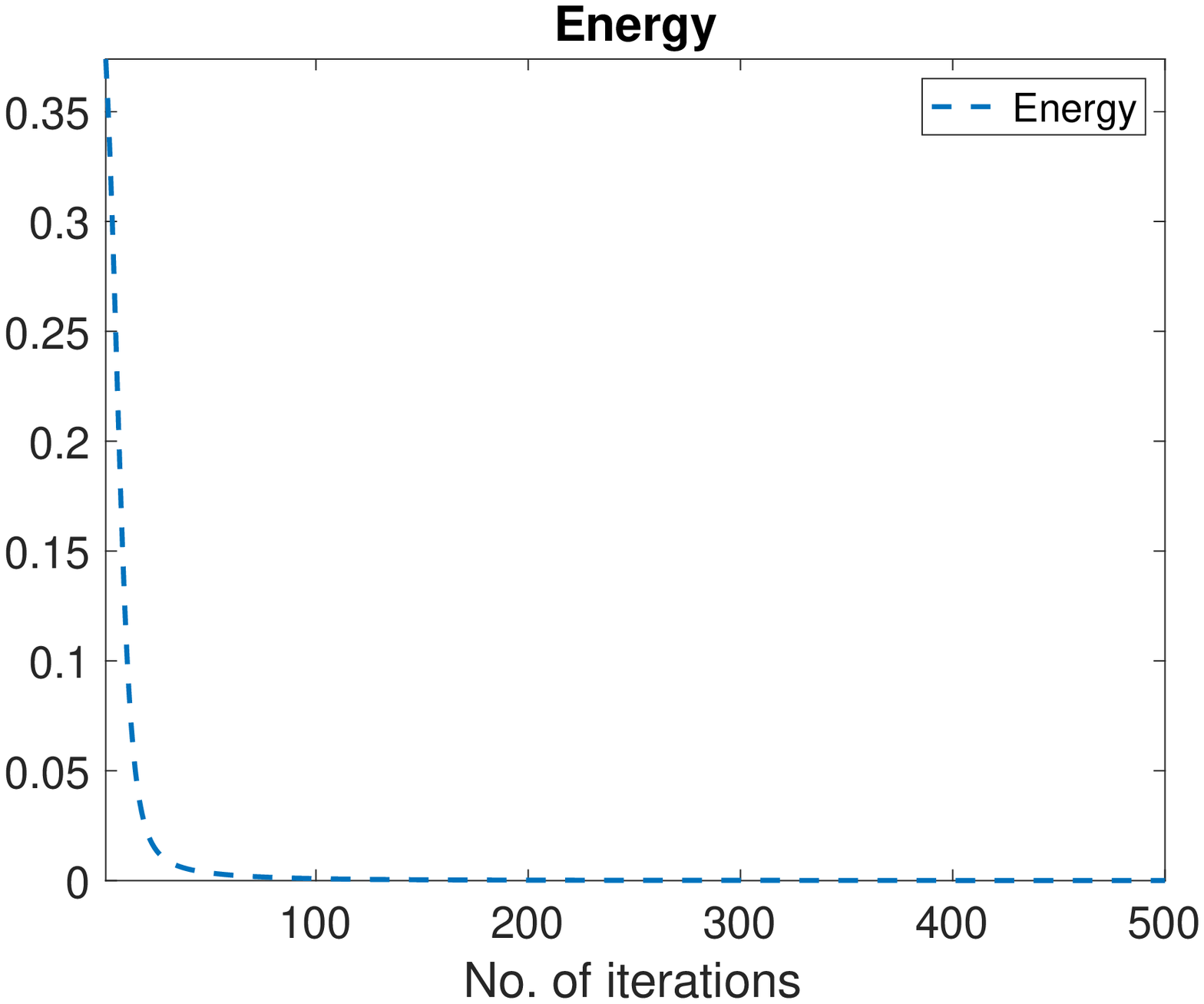}
\caption{Energy decrease of objective}
\label{subfig:results_gt2_2}
\end{subfigure}\\
\begin{subfigure}[b]{0.49\textwidth}
\centering
\includegraphics[width=\textwidth]{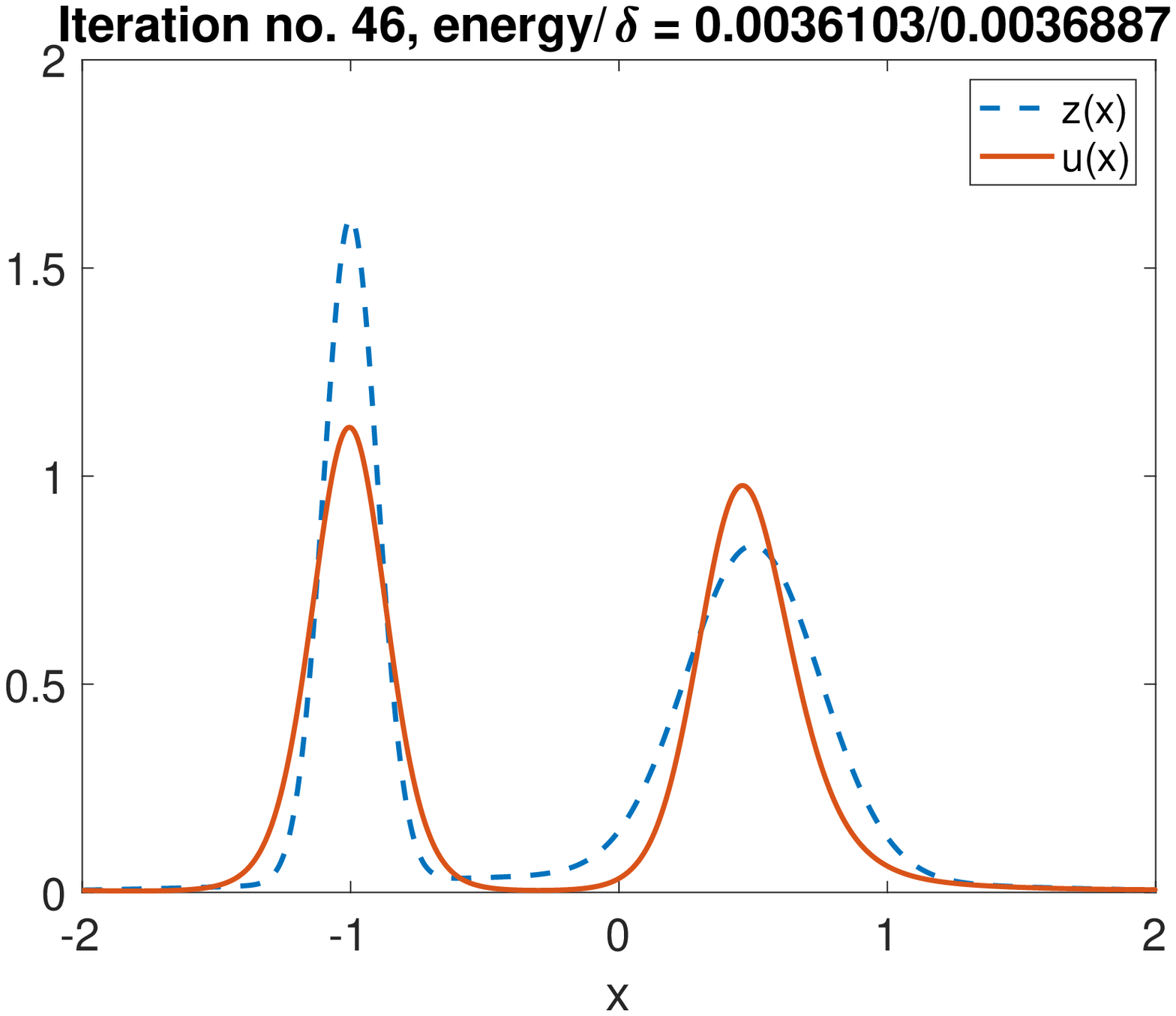}
\caption{$u^{k^\ast}$}
\label{subfig:results_gt2_3}
\end{subfigure}
\begin{subfigure}[b]{0.49\textwidth}
\centering
\includegraphics[width=\textwidth]{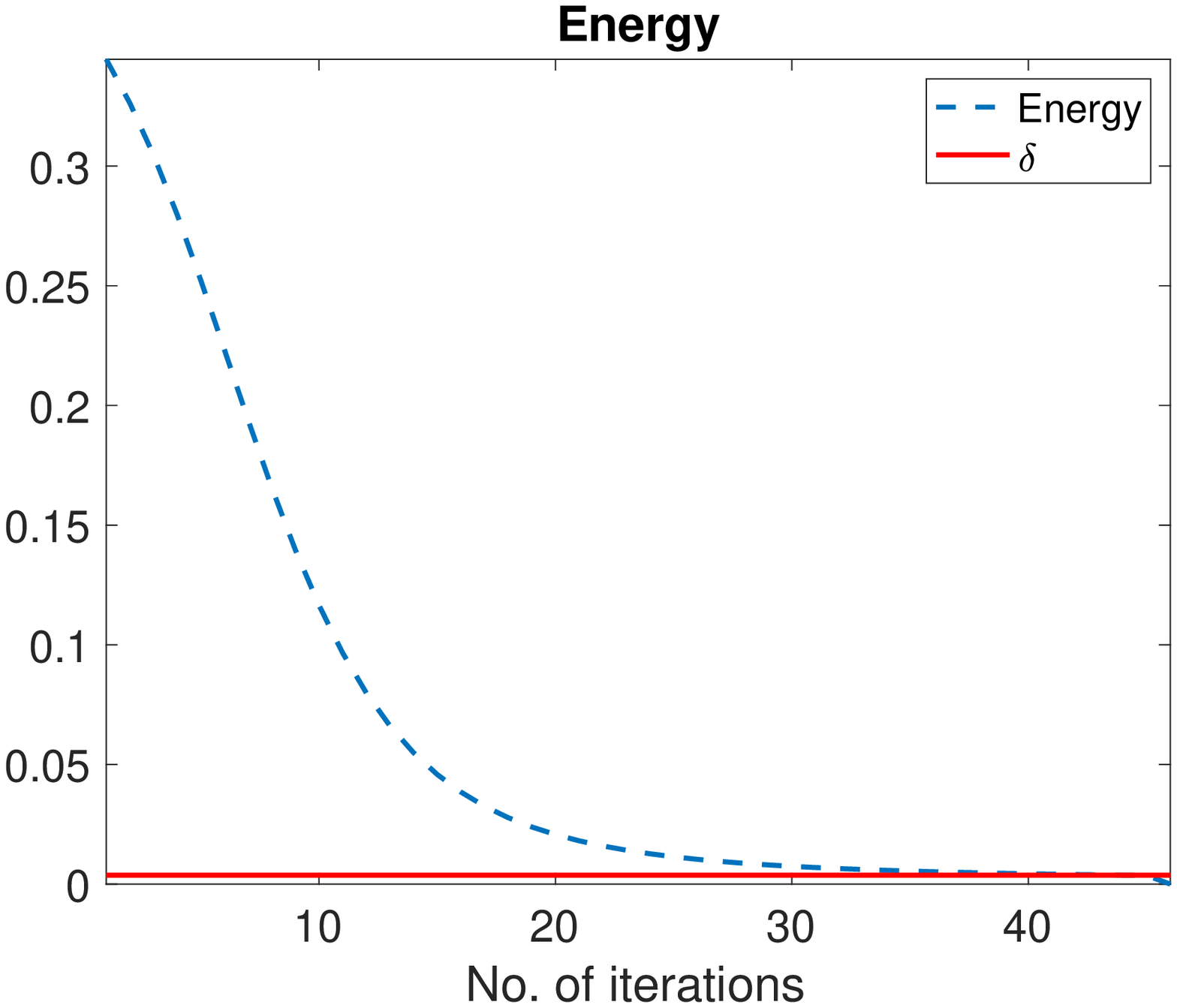}
\caption{Energy decrease of objective}
\label{subfig:results_gt2_4}
\end{subfigure}
\caption{Results of algorithm \eqref{method0b} for data measurements of the form \eqref{eq:measurements} based on the underlying function $z_2$ as visualised in Figure \ref{subfig:sampling_gt2}. Figure \ref{subfig:results_gt2_1} shows the ground truth function $z_2$ and $u^{k^\ast}$ under the assumption of $\sigma = 0$ in \eqref{eq:measurements} , whereas Figure \ref{subfig:results_gt2_2} visualises the monotonic energy decrease over the course of the iterations. Figure \ref{subfig:results_gt2_3} and Figure \ref{subfig:results_gt2_4} show the same results under the assumption of $\sigma = 1/500$ in \eqref{eq:measurements}.}
\label{fig:sampling_results2}
\end{figure}

Note that  $\tilde{z}$ does not satisfy the source condition \eqref{sc}, which is why we design a second function
\begin{align*}
z(x) := c \, \text{Re}\left( \frac{1}{\sqrt{2}} \sum_{j = 1}^n \left( \int_{\Omega} \tilde{z}(t) \, e^{-i t  \xi_j} \, dt \right) e^{i x \xi_j} \right) \, ,
\end{align*}
where $c$ is chosen to ensure $\int_{-10}^{10} z(x) \, dx = 1$, which by construction satisfies \eqref{sc}. We design two functions $z_1$ and $z_2$; $z_1$ is defined as $z_1 := z$ for $\tilde z$ with means $\mu_1 = 0, \mu_2 = -1$ and $\mu_3 = 1/2$, standard deviations $\sigma_1 = 1, \sigma_2 = 1/10$ and $\sigma_3 = 1/4$, and coefficients $c_1 = 1/10, c_2 = 3/5$ and $c_3 = 3/10$. The function $z_2 := \tilde z$ has the same means and standard deviations as $\tilde z$ in the previous example, but coefficients $c_1 = 1/10, c_2 = 2/5$ and $c_3 = 1/2$ instead. Both functions are visualised in Figure \ref{fig:sampling_ground_truth}. Subsequently we create data samples via \eqref{eq:measurements} with noise levels $\sigma = 0$ and $\sigma = 1/500$ for $z_1$, respectively $\sigma = 0$ and $\sigma = 1/50$ for $z_2$. 

\begin{figure}[t]
\begin{subfigure}[b]{0.49\textwidth}
\centering
\includegraphics[width=\textwidth]{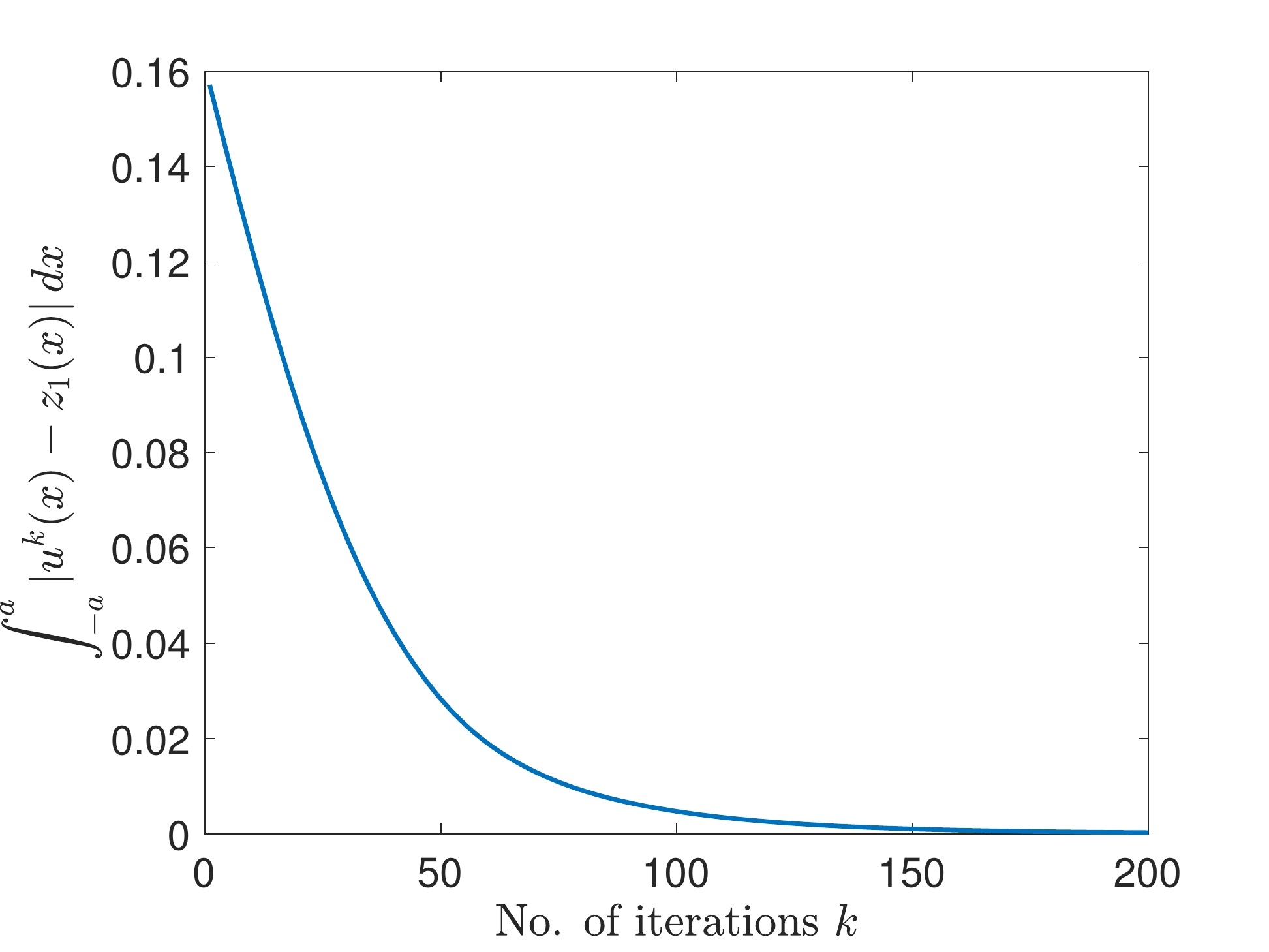}
\caption{}
\label{subfig:l1_difference_z1}
\end{subfigure}
\begin{subfigure}[b]{0.49\textwidth}
\centering
\includegraphics[width=\textwidth]{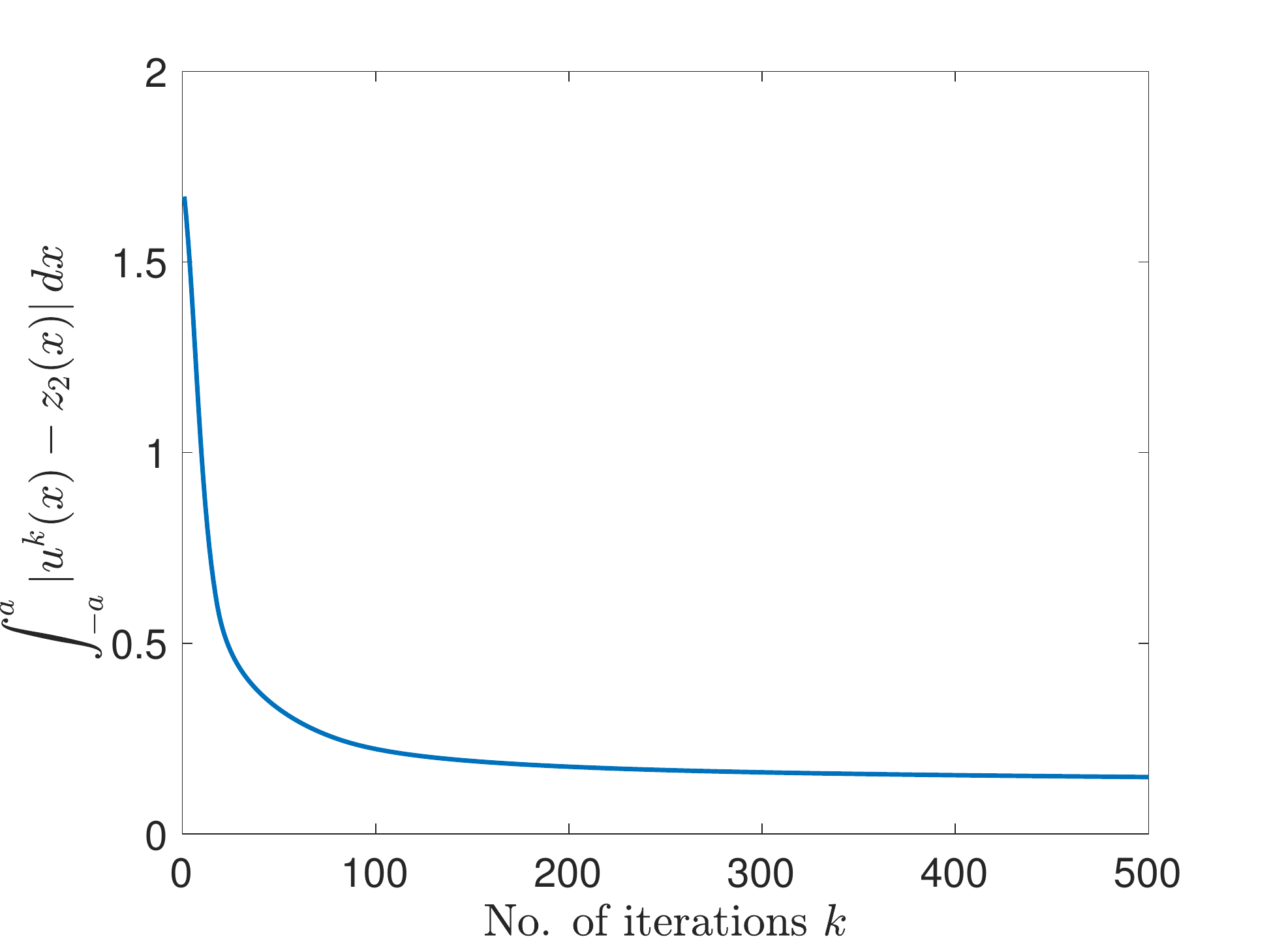}
\caption{}
\label{subfig:l1_difference_z2}
\end{subfigure}
\caption{The $L^1$-norm of the difference of the iterates $u^k$ and $z_1$ (Fig. \ref{subfig:l1_difference_z1}) and $z_2$, respectively (Fig. \ref{subfig:l1_difference_z2}), for the case of exact data ($\sigma = 0$). }
\label{fig:sampling_l1norm}
\end{figure}

In the following we run the entropic projection method \eqref{entropic} for $j = 1$, $\tau = 9/(10 \sqrt{2 \pi})$ and with the initial function
\begin{align*}
u_0(x) := \frac{1}{20} \begin{cases} 1 & x \in [-10, 10] \\ 0 & \text{else} \end{cases} \, ,
\end{align*}
either until the discrepancy principle \eqref{aposteriori} is violated (for $\tau = 1$) or until we reach a certain maximum number of iterations. We first investigate the algorithm for the function $z_1$ as seen in Figure \ref{subfig:sampling_gt1}, for perfect data ($\sigma = 0$) and for noisy data $\sigma = 1/500$. For perfect data we run the algorithm for 201 iterations and observe that we are converging towards $z_1$ as can be seen in Figure \ref{subfig:results_gt1_1} as well as in Figure \ref{subfig:l1_difference_z1}, which is a numerical confirmation of Proposition \ref{pro:exact_data}. For the non-trivial noise-level $\sigma = 1/500$ the algorithm stops after 65 iterations according to the discrepancy principle (Figure \ref{subfig:results_gt1_3}).

To conclude, we run the same numerical experiments for $z_2$ as shown in Figure \ref{subfig:sampling_gt2}. As we mentioned earlier, \eqref{sc} is violated and even for perfect data (i.e. $\sigma = 0$) we cannot expect the results of Proposition \ref{pro:exact_data} to hold true. It can be seen in Figure \ref{subfig:results_gt2_1} that $u^k$ does not seem to converge towards $z_2$ despite a decrease of the objective to values in the order of $10^{-5}$. In fact, if we compare the $L^1$-norm of the difference $u^k - z_2$, we also see in Figure \ref{subfig:l1_difference_z2} that $u^k$ does not seem to converge towards $z_2$. For noisy data with $\sigma = 1/500$ the discrepancy principle is violated after 46 iterations, with its result being visualised in Figure \ref{subfig:results_gt2_3}.

\subsection{Initial Densities for Stochastic Differential Equations}

An interesting problem in several applications, e.g. in data assimilation scenarios (cf. e.g \cite{Hairer}), is the reconstruction of the initial density for a system evolving via stochastic differential equations with drift $b$ and volatility $a$. The density evolves via the Fokker-Planck equation (cf. \cite{gardiner})
\begin{equation}
	\partial_t \rho(x,t) + \nabla \dot (\rho(x,t) b(x,t)) =  \frac{1}2 \Delta( a(x,t)^2 \rho(x,t))  
\end{equation}
in $\Omega \times (0,T)$ with no-flux boundary conditions. Under appropriate smoothness conditions on $a$ and $b$  as well as positivity of $a$ is is well-known that the Fokker-Planck equation has a unique nonnegative solution $
\rho \in C(0,T;L^1(\Omega))$ for nonnegative initial values $u \in L^1(\Omega)$ such that $\int_\Omega u \ln u ~dx < \infty$. In problems related to reconstructing $u$ it is hence rather natural to use methods penalizing its entropy. 

The forward operator $A$ maps the initial density to indirect measurements of the density $\rho$ over time, e.g. moments or local integrals. Parametrizing the measurements by values $\sigma$ in a bounded set $\sigma$ we obtain 
\begin{equation}
	A: L^1(\Omega) \rightarrow L^2((0,T) \times \Sigma ) , \qquad u \mapsto \int_\Omega k(\sigma,y) \rho(y,t)~dy .
\end{equation}
It is well-known that Fokker-Planck equations satisfy an $L^1$-contractivity property on the domain of the entropy functional (cf. \cite{Karlsen}) i.e. for $\rho_i$ denoting the solution with initial value 
\begin{equation}
	\Vert \rho_1(t) - \rho_2(t) \Vert_{L^1} \leq \Vert u_1 - u_2 \Vert
\end{equation}
for almost all $t \in (0,T)$. Thus, the map $u \mapsto \rho$ is Lipschitz continuous with unit modulus when considered as a map into $L^\infty(0,T;L^1(\Omega))$ on the domain of the entropy. Hence, if $k \in L^\infty(\Sigma \times \Omega)$ we can easily verify that the operator $A$ satisfies \eqref{condition}.

We finally mention that in the case of stationary coefficients $a$ and $b$, the Fokker-Planck equation has a unique stationary solution $\rho_\infty$ among nonnegative functions with unit mass (cf. \cite{droniou}), to which it converges with exponential speed in the relative entropy (cf. \cite{Arnold1,Arnold}), i.e. the Bregman distance related to the entropy functional. Hence, it is natural to use $\rho_\infty$ as an initial value for the reconstruction of $u$, since we may expect them to be close in particular in the relative entropy.

\section{Conclusions and remarks}
 
In this study we have investigated a multiplicative  entropic type method for ill-posed equations, which preserves nonnegativity of the iterates.  Historically, the underlying strategy has spreading roots in the inverse problems literature: Landweber iterates, surrogate functionals and linearized Bregman, to quote  a few approaches. In parallel, this has been treated in different contexts in finite dimensional optimization, e.g., as a mirror descent or as a steepest descent (linearized proximal)  algorithm with generalized distances, or in machine learning - as an exponentiated gradient descent algorithm for online prediction via linear models.

The closed form algorithm is shown to converge weakly in $L^1$ to a solution of the ill-posed problem and convergence rates are obtained by means of the Kullback-Leibler (KL) distance. 
All the results are quite naturally established when imposing "mean one" restriction to the unknown, while the case without restrictions relies on a norm combined with KL distance based  Lipschitz condition, in which case operators satisfying it remain to be found. 


Methods  of this type involving other interesting fidelity terms, nonlinear  operators and  eventually stochastic versions and line search  
strategies might be considered in more detail for future research.

%
%
%
%
%
%

\section{Acknowledgements}
Martin Burger acknowledges support from European Union’s Horizon 2020 research and innovation programme under the Marie Sk lodowska-Curie grant agreement No 777826 (NoMADS).
 Martin Benning  acknowledges support from the Leverhulme Trust Early Career Fellowship ECF-2016-611 'Learning from mistakes: a supervised feedback-loop for imaging applications'.


\end{document}